\documentstyle[12pt,amscd,amssymb]{amsart} 
        \newcommand{\rank}{\operatorname{rank}}

        \newcommand{\pardeg}{\operatorname{pardeg}}
        \newcommand{\End}{\operatorname{End}}
        \newcommand{\ParHom}{\operatorname{ParHom}}
        \newcommand{\ParEnd}{\operatorname{ParEnd}}
        \newcommand{\ParAut}{\operatorname{ParAut}}
        \newcommand{\Ext}{\operatorname{Ext}}
        \newcommand{\parhom}{\operatorname{{\frak P}{\frak a}{\frak r}{\frak H}{\frak o}{\frak m}}}
        \newcommand{\parend}{\operatorname{{\frak P}{\frak a}{\frak r}{\frak E}{\frak n}{\frak d}}}
        \newcommand{\gr}{\operatorname{gr}}
        \renewcommand{\P}{\operatorname{P}}
        
        \newcommand{\Res}{\operatornamewithlimits{Res}}
        \renewcommand{\aa}{\alpha}
        
        \newcommand{\delbar}{\overline{\partial}}
        \newcommand{\Tr}{\operatorname{Tr}}
        \newcommand{\Ker}{\operatorname{Ker}}
        \newcommand{\Aut}{\operatorname{Aut}}

	\newtheorem{defn}{Definition}[section]
	\newtheorem{lem}[defn]{Lemma}
	\newtheorem{thm}[defn]{Theorem}
	\newtheorem{prop}[defn]{Proposition}
	\newtheorem{cor}[defn]{Corollary}
	\newtheorem{clm}[defn]{Claim}

\begin{document}
\title{Moduli Spaces of Parabolic Higgs Bundles and Parabolic K(D) Pairs over Smooth Curves: I}
\author{Hans U. Boden and K\^oji Yokogawa} 
\address{Max-Planck-Institut f\"ur Mathematik \\ Gottfried-Claren-Strasse 26 \\ 53225 Bonn \\ Germany}
\email{hboden@@mpim-bonn.mpg.de}
\email{yokogawa@@mpim-bonn.mpg.de}
\begin{abstract}
This paper concerns the moduli spaces of rank two 
parabolic Higgs bundles and parabolic $K(D)$ pairs over a smooth curve.
Precisely which parabolic bundles occur 
in stable $K(D)$ pairs and stable Higgs bundles is determined.
Using Morse theory, the moduli space of parabolic Higgs bundles 
is shown to be a non-compact, connected, simply connected manifold, and
a computation of its Poincar\'e polynomial is given. 
\end{abstract}
 
\maketitle

\section{Introduction}
Let $C$ be a compact curve. The correspondence between unitary representations of $\pi_1(C)$ 
and  semistable bundles over $C$ of degree zero  
\cite{narasimhan-seshadri} was extended to
non-compact curves $C_0$ by Mehta and Seshadri in \cite{mehta-seshadri}.
If $C_0$ has
compactification $C,$
they
prove that semistable parabolic bundles over $C$  of parabolic degree zero 
correspond to unitary representations of $\pi_1(C_0)$ with fixed holonomy around
$p \in C \setminus C_0.$ 
{Generalizing in a different direction, Hitchin and Donaldson
\cite{hitchin, donaldson3}
proved that
representations of $\pi_1(C)$ correspond to
semistable Higgs bundles over $C$ of degree zero.}\footnote{The non-abelian Hodge theorem, a further generalization of this,
holds for arbitrary compact K\"ahler manifolds
\cite{donaldson1, donaldson2, uhlenbeck-yau, corlette,simp1,simp4}.}
A Higgs bundle includes the additional information of a Higgs field,
which is a holomorphic map $\Phi : E \rightarrow E \otimes K,$ where $K$ denotes the canonical bundle.
 
In the case of a parabolic bundle,
the Higgs field is permitted to have poles
of order one at the compactification points.
Requiring these residues
to be either {\it parabolic} or {\it nilpotent},
one obtains
two moduli spaces: ${\cal P}_\aa,$ 
the moduli space of parabolic $K(D)$ pairs,
and ${\cal N}_\aa,$ the moduli space of parabolic Higgs bundles.
The subscript $\alpha$ refers to a particular choice of weights.
In \cite{yoko1}, ${\cal P}_\aa$ is constructed using Geometric Invariant Theory and
is proved to be a normal, quasi-projective variety.
In \cite{konno},  ${\cal N}_\aa$ is constructed as a hyperk\"ahler quotient using gauge theory.

Simpson's factorization theorem
states that for $X$ a projective  algebraic variety,
any $SL(2, {\Bbb C})$ representation of $\pi_1(X)$ with
Zariski dense image is either rigid or factors through an algebraic map
from $X$ to an orbicurve \cite{simp3}.
Because orbicurve representations can be interpreted
as stable parabolic Higgs bundles \cite{scheinost},
it is important to understand these moduli spaces,
which is the subject of our study here.

Given a rank two parabolic bundle,
we first establish algebraic conditions
for the existence of
a field making it stable as either a $K(D)$ pair or
a Higgs bundle. One could use this to describe
both moduli spaces,
which we do for one particular case, but this approach appears too complicated 
to work in general. 

For that reason, we shift gears and study the topological properties of the moduli space of 
parabolic Higgs bundles, using the approach of Hitchin \cite{hitchin}.
There is a circle action on  ${\cal N}_\aa$ preserving 
its complex and symplectic structure and
the associated moment map is a Morse function in the sense of Bott.
We prove that ${\cal N}_\aa$ is a {non-compact}\footnote{The exception to
this case is studied in \S 2.3}, connected, simply connected manifold
and compute its Betti numbers, which turn out to be independent of 
the weights $\aa.$
This is surprising because it is not true for non-Higgs bundles:
the Betti numbers of the moduli space ${\cal M}_\aa$
of parabolic bundles do depend
in an essential way on $\aa$ (cf. \cite{bh}).
In the sequel, we plan to extend these results to higher rank bundles. 

The paper is organized as follows.
In \S 2.1 we define parabolic bundles with auxiliary fields and introduce
the three moduli spaces ${\cal M}_\alpha, \, {\cal N}_\alpha,$ and ${\cal P}_\alpha.$
Tensor products, duals, and the
Serre duality theorem for parabolic bundles are given in
\S 2.2.
In \S 3.1, we establish the algebraic conditions mentioned above, and in \S 3.2,
we use these 
to characterize ${\cal P}_\alpha$ and ${\cal N}_\alpha$
in the case of ${\Bbb P}^1$ with three parabolic points.
Turning our attention to ${\cal N}_\alpha$ in \S 4,
we describe its construction in \S 4.1 
as a hyperk\"ahler quotient, following \cite{konno}.
In \S 4.2 we define the Morse function 
on ${\cal N}_\alpha^0$ and then prove our main results 
about the topology of ${\cal N}_\alpha^0$ in \S\S 4.3 and 4.4.

Both authors are grateful to the Max-Planck-Institut f\"ur Mathematik for providing
a stimulating intellectual environment as well as financial support.
Warm thanks also to the VBAC Research Group of Europroj for travel funding, 
and to O. Garcia-Prada, L. G\"ottsche, N. Hitchin, Y. Hu, and D. Huybrechts 
for their advice.

After submitting this paper, we learned that Nasatyr and
Steer have obtained similar results
studying orbifold Higgs bundles \cite{nasatyr-steer}. 

\section{Definitions and Preliminary Results}
\subsection{Three moduli spaces}

Let $X$ be a smooth curve of genus $g$ with $n$ marked points
in the reduced divisor $D=p_1+\cdots + p_n$ and $E$ a holomorphic bundle over $X.$ 
\begin{defn} A parabolic structure on $E$ consists of
weighted flags 
\begin{eqnarray*}
&E_{p}=F_1(p) \supset \cdots \supset F_{s_p}(p) \supset 0&\\
&0\leq \alpha_1(p) < \cdots <\alpha_{s_p}(p) < 1&
\end{eqnarray*}
over each $p \in D.$
A holomorphic map $\phi:E^1 \longrightarrow E^2$ between parabolic bundles is called 
parabolic if
$\alpha^1_i(p) > \alpha^2_j(p)$ implies $\phi(F^1_i(p)) \subset F^2_{j+1}(p)$
for all $p \in D.$
We call $\phi$ strongly parabolic if 
$\alpha^1_i(p) \geq \alpha^2_j(p)$ implies $\phi(F^1_i(p)) \subset F^2_{j+1}(p)$
for all $p \in D.$
\label{defn:pbundle}
\end{defn}
We use $E_*$ to denote the bundle together with a parabolic structure.
Also, we use $\ParHom(E^1_*,E^2_*)$ and $\ParHom(E^1_*,\widehat{E}^2_*)$ to
denote the sets of parabolic and strongly
parabolic morphisms from $E^1$ to $E^2,$ 
respectively.
(The decorative notation will become clear in \S 2.2.)
If $\alpha^1_i(p) \neq \alpha^2_j(p)$ for all $i,j$ and $p \in D,$
then a parabolic morphism is automatically strongly parabolic.
On the other hand, using the notation 
$\ParEnd(E_*) = \ParHom(E_*,E_*)$ and $\ParEnd^{\wedge}(E_*) = \ParHom(E_*,\widehat{E}_*),$ 
then strongly parabolic endomorphisms are nilpotent with respect to the flag data at each $p \in D.$

Let $K$ denote the canonical bundle of $X$ and give $E\otimes K(D)$ the obvious parabolic structure. 
\begin{defn}
A parabolic $K(D)$ pair is a pair $(E, \Phi)$ consisting of  a parabolic bundle $E$
and a parabolic map $\Phi:E \rightarrow E \otimes K(D).$ 
Such a pair is called a parabolic Higgs bundle if, in addition,
$\Phi$ is a strongly parabolic morphism.
\end{defn}

Viewing $\alpha$ as a vector-valued function on $D,$ we use it as an index 
to indicate the parabolic structure on $E_*.$
Let $m_i(p) =\dim(F_i(p))-\dim(F_{i+1}(p)),$
the multiplicity of $\alpha_i(p),$ 
and $f_p = \frac{1}{2}(r^2 - \sum_{i=1}^{s_p} (m_i(p))^2),$ the dimension of the 
associated flag variety.
Define the parabolic degree and slope of $E_*$ by
\begin{eqnarray*}
\pardeg E_* &=& \deg E + \sum_{p \in D} \sum_{i=1}^{s_p} m_i(p) \alpha_i(p),\\
\mu (E_*) &=& {{\pardeg E_*}\over{\rank E}}.
\end{eqnarray*}

If $L$ is a subbundle of $E,$ then $L$ inherits a parabolic structure
from $E$ by pullback.
We call the bundle $E_*$ {\it stable} 
({\it semistable}) if, for every proper
subbundle $L$ of $E,$ we have $\mu (L_*) < \mu (E_*)$ 
(respectively $\mu (L_*) \leq \mu (E_*)$).
Likewise, we will call a parabolic $K(D)$ pair $(E_*,\Phi)$ {\it stable} (or {\it semistable})
if the same inequalities hold on those proper subbundles $L$ of $E$ which are, in addition, $\Phi$-invariant. 

Denote by ${\cal M}_\alpha$ the moduli space of $\alpha$-semistable parabolic bundles, 
by ${\cal N}_\alpha$ the moduli space of $\alpha$-semistable parabolic Higgs bundles, and
by ${\cal P}_\alpha$ the moduli space of $\alpha$-semistable parabolic $K(D)$ pairs.
By \cite{mehta-seshadri}, ${\cal M}_\alpha$ is a normal, projective variety of dimension
$$\dim {\cal M}_\alpha = (g-1)r^2+1+\sum_{p \in D} f_p .$$
(If $g=0,$ this holds only when ${\cal M}_\alpha \neq \emptyset.$) 
Further, in \cite{yoko1,yoko2}, ${\cal P}_\alpha$ is shown to be a normal, 
quasi-projective variety of dimension
$$\dim {\cal P}_\alpha = (2g-2+n)r^2+1$$
which contains ${\cal N}_\alpha$ 
as a closed subvariety of ${\cal P}_\alpha$ 
of dimension
$$\dim {\cal N}_\alpha = 2(g-1)r^2+2+2\sum_{p \in D} f_p.$$

For generic $\alpha,$ a bundle (or pair) 
is $\alpha$-semistable $\Leftrightarrow$ it is $\alpha$-stable.
In these cases, the moduli spaces ${\cal M}_\alpha, {\cal N}_\alpha$ and ${\cal P}_\alpha$ are smooth and
can be described topologically as certain quotients of the gauge group
${\cal G}^{\Bbb C}=\ParAut(E_*).$ 
The same is true for ${\cal M}^0_\alpha, {\cal N}^0_\alpha$ and ${\cal P}^0_\alpha,$
the moduli spaces with fixed determinant and trace-free $\Phi.$
In this way, it is shown in \cite{konno} that ${\cal N}^0_\alpha$ is, for generic $\alpha,$ a smooth,
hyperk\"{a}hler manifold of complex dimension
$$\dim {\cal N}^0_\alpha = 2(g-1)(r^2-1)+2\sum_{p \in D} f_p .$$

\subsection{Parabolic sheaves and Serre duality}
Some of the material in this section is a summary of results in \cite{yoko2}

Suppose now that $E$ is a locally free sheaf on $X$ and $D=p_1+\cdots + p_n$ is a reduced divisor.
\begin{defn}
A parabolic structure on $E$ consists of a weighted filtration of the form 
\begin{eqnarray*} 
&E=E_0 = E_{\alpha_1} \supset \cdots \supset E_{\alpha_l} \supset E_{\alpha_{l+1}} = E(-D),&\\
&0=\alpha_0\leq \alpha_1 < \cdots < \alpha_l < \alpha_{l+1}=1.&
\end{eqnarray*}
We can define $E_x$ for $x \in [0,1]$ by setting $E_x = E_{\alpha_i}$
if $\alpha_{i-1} < x \leq \alpha_i,$ and then extend to $x \in {\Bbb R}$ by setting $E_{x+1} =E_{x}(-D).$ 
We call the resulting filtered sheaf $E_*$ a parabolic sheaf.

We define the coparabolic sheaf $\widehat{E}_*,$ by 
$$\widehat{E}_x = \left\{\begin{array} {ll} E_x & \hbox{ if } x \neq \alpha_i\\
                                        E_{\alpha_{i+1}} & \hbox{ if } x = \alpha_i.
\end{array} \right.$$
A morphism of parabolic sheaves $\phi:E^1_*\rightarrow E^2_*$ is
a called parabolic if $\phi(E^1_x) \subseteq E^2_x$
and strongly parabolic if $\phi(E^1_x) \subseteq \widehat{E}^2_x$
for all $x \in {\Bbb R}$.
\end{defn}
\begin{figure}[b]
\begin{picture}(-80,100)(250,0)
        \put(-20,50){$E_*$} \put(0,10){\line(1,0){160}} \put(20,0){\line(0,1){90}} \put(120,0){\line(0,1){90}} \put(22,0){$0$} \put(122,0){$1$} \put(30,85){$E$} \put(0,80){\line(1,0){38}} \put(40,80){\circle*{4}} \put(40,0){$\aa_1$} \multiput(40,10)(0,4){12}{\line(0,1){2}} \multiput(40,78)(0,-4){4}{\line(0,-1){2}}\put(60,65){$E_{\aa_2}$} \put(42,60){\line(1,0){28}} \put(40,60){\circle{4}} \put(70,60){\circle*{4}} \put(70,0){$\aa_2$} \multiput(70,10)(0,4){7}{\line(0,1){2}} \multiput(70,58)(0,-4){4}{\line(0,-1){2}}\put(72,40){\line(1,0){28}} \put(90,45){$E_{\aa_3}$} \put(70,40){\circle{4}} \put(100,40){\circle*{4}} \put(100,0){$\aa_3$} \multiput(100,10)(0,4){5}{\line(0,1){2}} \multiput(100,38)(0,-4){2}{\line(0,-1){2}} \put(102,30){\line(1,0){38}} \put(100,30){\circle{4}} \put(140,30){\circle*{4}} \put(125,35){$E(-D)$} \put(140,0){$1+\aa_1$} \multiput(140,10)(0,4){2}{\line(0,1){2}} \multiput(140,28)(0,-4){2}{\line(0,-1){2}} \put(142,20){\line(1,0){18}} \put(140,20){\circle{4}}
        \put(180,50){$\Rightarrow$} \put(200,50){$\widehat{E}_*$} \put(220,10){\line(1,0){160}} \put(240,0){\line(0,1){90}} \put(340,0){\line(0,1){90}} \put(242,0){$0$} \put(342,0){$1$} \put(250,85){$E$} \put(220,80){\line(1,0){38}} \put(260,80){\circle{4}} \put(260,0){$\aa_1$} \multiput(260,10)(0,4){17}{\line(0,1){2}} \put(280,65){$E_{\aa_2}$} \put(260,60){\line(1,0){28}} \put(260,60){\circle*{4}} \put(290,60){\circle{4}} \put(290,0){$\aa_2$} \multiput(290,10)(0,4){12}{\line(0,1){2}} \put(290,40){\line(1,0){28}} \put(310,45){$E_{\aa_3}$} \put(290,40){\circle*{4}} \put(320,40){\circle{4}} \put(320,0){$\aa_3$} \multiput(320,10)(0,4){7}{\line(0,1){2}} \put(320,30){\line(1,0){38}} \put(320,30){\circle*{4}} \put(360,30){\circle{4}} \put(345,35){$E(-D)$} \put(360,0){$1+\aa_1$} \multiput(360,10)(0,4){5}{\line(0,1){2}} \put(360,20){\line(1,0){18}} \put(360,20){\circle*{4}}
\end{picture}
\caption{The simple relationship between $E_*$ and $\widehat{E}_*.$}
\end{figure}
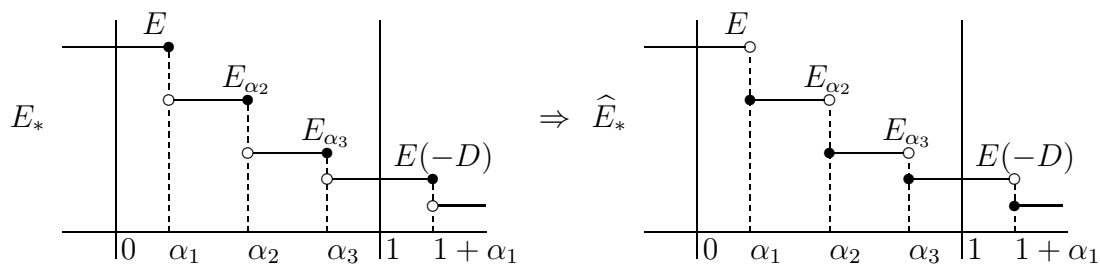

We shall denote by
$\parhom(E^1_*, E^2_*)$ and $\parhom(E^1_*, \widehat{E}^2_*)$ 
the sheaves of parabolic and strongly parabolic morphisms,
and by $\ParHom(E^1_*,E^2_*)$ and $\ParHom(E^1_*, \widehat{E}^2_*)$
their global sections. 
We now show that there is an equivalence of the categories of
parabolic bundles on $X$ and parabolic sheaves on $X.$

Given a parabolic bundle $E$ with flags and weights as in Definition \ref{defn:pbundle},
we define the  filtered sheaf $E_*$ 
following Simpson \cite{simp5}. 
For $p \in D$ and $\alpha_{i-1}(p) < x \leq\alpha_i(p)$, set 
\begin{eqnarray*}
E_x^p &=& \ker(E\rightarrow E_{p}/F_i(p)),\\ 
E_x &=& \bigcap_{p \in D} E^p_x.
\end{eqnarray*}
Now extend to all $x$ by  $E_{x+1} = E_x(-D).$

Conversely, given a parabolic sheaf $E_*,$ the quotient
$E/E_1$ is a skyscraper sheaf with support on $D$ and, for each $p \in D,$
we get weighted flags in $E_p$ by intersecting with the filtration at $p.$ 
To be precise, let 
$\alpha_1(p), \ldots, \alpha_{s_p}(p)$ be the subset of weights
such that 
\begin{equation}
\alpha_{i-1}(p) < x \leq \alpha_{i}(p) \Leftrightarrow (E_x/E_1)_p = (E_{\alpha_{i}(p)}/E_1)_p.
\label{eqn:1}
\end{equation}
Setting $F_i(p) = (E_{\alpha_{i}(p)}/E_1)_p,$
we obtain a parabolic bundle in the sense of Definition \ref{defn:pbundle}.

Suppose now $E^1$ and $E^2$ are parabolic bundles and $\phi\in\ParHom(E^1,E^2).$ 
We want to show that $\phi$ induces a morphism of the parabolic sheaves.
So, suppose $\alpha^1_{i-1}(p) < x \leq\alpha^1_i(p)$ and
$\alpha^2_{j-1}(p) < x \leq\alpha^2_j(p).$
Since $\alpha^1_i(p) > \alpha^2_{j-1}(p),\; \phi(F^1_i(p)) \subset F^2_j(p)$ 
and we see that $\phi$ maps $\ker(E^1\rightarrow E^1_{p}/F^1_i(p))$
to $\ker(E^2\rightarrow E^2_{p}/F^2_i(p))$ for all $p \in D,$
from which it follows that $\phi$ induces a map $\phi:E^1_x \rightarrow E^2_x.$

Suppose conversely that $E^1_*$ and $E^2_*$ are parabolic sheaves, $\phi \in \ParHom(E^1_*,E^2_*)$
and $\alpha^1_i(p) > \alpha^2_j(p).$ Set $x=\alpha^1_i(p)$ and
$y=\alpha^2_{j+1}(p)$ for notational convenience. 
Then $\phi(E^1_x) \subset E^2_x.$ Since $x > \alpha^2_j(p),$ it follows from (\ref{eqn:1}) that
$(E^2_x/E^2_1)_p \subset (E^2_{y}/E^2_1)_p$ and hence
$\phi(F^1_i(p)) \subset F^2_{j+1}(p).$ 

It is not hard to see the same correspondence for strongly parabolic morphisms. 
Thus, we have an equivalence of the categories of parabolic bundles and
parabolic sheaves. We use the definitions interchangeably and
denote by $E_*$ a parabolic bundle or sheaf,
reserving $E = E_0$ for the underlying holomorphic bundle.
\medskip

For the convenience of readers, we briefly summarize the results in \cite{yoko2}
dealing with exact sequences and  tensor products of parabolic sheaves.
This is necessary for the statement of Serre duality for parabolic bundles,
which is a tool we use throughout the paper.

The category of parabolic sheaves ${\cal P}$ is not abelian, but is contained in
an abelian category $\widetilde{{\cal P}}$ as a full subcategory. 
Objects in $\widetilde{{\cal P}}$ are also written by $E_*$ and a morphism 
$f:E_*^1\rightarrow E_*^2$ is a family of morphisms $f_x:E_x^1\rightarrow E_x^2.$ 
A coparabolic sheaf $\widehat{E}_*$ is realized 
in $\widetilde{{\cal P}}$. The set 
$\ParHom(E^1_*, \widehat{E}^2_*)$ is just the set of morphisms in 
$\widetilde{{\cal P}}.$
In $\widetilde{{\cal P}},$  a sequence 
\begin{equation}
	0 \longrightarrow L_*
	\longrightarrow E_*
	\longrightarrow M_*
	\longrightarrow 0
\label{eqn:a}
\end{equation}
is exact if and only if the induced sequence at $x$ is exact for all $x\in {\Bbb R}.$ 

\medskip
\noindent
{\it Remark. \,} 
If the sequence (\ref{eqn:a}) is exact, then so is the sequence obtained by tensoring 
(\ref{eqn:a}) with any parabolic bundle  (cf. Proposition 3.3 of \cite{yoko2}) and 
$$\pardeg E_* = \pardeg L_* +\pardeg M_*.$$

We can define dual parabolic sheaves $E_*^\vee$, parabolic tensor products $L_*\otimes M_*$, 
Hom-parabolic sheaves $\parhom(L_*, M_*)_*$, 
and cohomology groups $\Ext^i(L_*,M_*)$. 
Clearly, 
$$\pardeg (L_*\otimes M_*) = \rank(M)\pardeg L_* + \rank(L)\pardeg M_*.$$ 
In addition, we have  
\begin{eqnarray*}
\Ext^0(L_*,M_*) &=& H^0(L_*^\vee \otimes M_*) = H^0(\parhom(L_*, M_*)) = \ParHom(L_*, M_*),\\
\Ext^1(L_*,M_*) &=& H^1(L_*^\vee \otimes M_*) = H^1(\parhom(L_*,M_*)). 
\end{eqnarray*}
We can identify $\Ext^1(M_*,L_*)$ with the set of equivalence classes of exact sequences of type (\ref{eqn:a}). 

The Serre duality theorem is generalized as follows (see Proposition 3.7 of \cite{yoko2}).  

\begin{prop}\label{serre-dual-thm}
For parabolic sheaves $L_*$ and $M_*$, there is a natural isomorphism
$$	\theta^i:H^i(L_*^\vee \otimes M_* \otimes K(D))
        \stackrel{\simeq}{\longrightarrow} 
	H^{1-i}(M_*^\vee \otimes \widehat{L}_*)^\vee.$$
\end{prop}

Given $E_*$ and $\beta \in {\Bbb R}^n,$
define 
$E_*[\beta]_*,$ 
the parabolic sheaf $E_*$ shifted by $\beta,$ 
by 
$$E_*[\beta]_x = \bigcap_{i} E^{p_i}_{x+\beta_i}.$$
{\it Example. \;} {\it The Picard group of parabolic line bundles.}\\
A holomorphic bundle $E$ is regarded as a parabolic bundle with 
the trivial parabolic structure $E_p\supset 0, \aa_1(p)=0$ at each $p \in D.$
We call this the special structure on $E$. 
Note that every parabolic line bundle $L_*$ is gotten by shifting 
the special structure on the underlying bundle $L,$ i.e., there is
a unique $\beta \in [0,1)^n$ with
$L_*=L[\beta]_*$
Viewing ${\cal O}_X$ as a parabolic bundle with the special structure,
then it is not difficult to verify that
\begin{equation}\label{eqn:b}
E_*[\beta]_* = E_* \otimes {\cal O}_X[\beta]_*
\end{equation}
Let $e_i$ denote the standard basis vector in ${\Bbb R}^n.$
From (\ref{eqn:b}) we have
\begin{eqnarray*}
E^1_*[\beta^1]_* \otimes E^2_*[\beta^2]_* &=& E^1_*\otimes E^2_*[\beta^1 + \beta^2]_*, \\
E_*[\beta]_*^\vee &=& E^\vee_*[-\beta]_*, \\
E_*[e_i]_* &=& E_*\otimes{\cal O}_X(-p_i). 
\end{eqnarray*}
These three formulas determine 
the Picard group of parabolic line bundles on $X$.

\medskip
\noindent
{\it Remark. \,}
For any parabolic line bundle $L_*$,
the stability (or semistability) of $E_*\otimes L_*$
is equivalent to that of $E_*.$
Similarly, the stability (or semistability) of 
$(E_*\otimes L_*,\Phi\otimes 1)$ 
is equivalent to that of $(E_*, \Phi).$

\medskip

In particular, apply this to the case  
of a rank two parabolic bundle $E_*$ with full flags at each $p_i$
and weights $0 \leq \alpha_1(p_i) < \alpha_2(p_i) < 1.$
Using equation (\ref{eqn:b}) with $\beta_i =\frac{1}{2}(\aa_1(p_i)+\aa_2(p_i)-1),$
notice that $E_*[\beta]_*$ has weights
$0 < a_1(p_i) < 1-a_1(p_i) < 1$ at $p_i,$
where $a_1(p_i) = \frac{1}{2}(\aa_1(p_i)-\aa_2(p_i)+1).$

\section{An Algebraic Description of the Moduli Spaces in Rank Two}
\subsection{Criteria for the existence of stabilizing fields}
In this section, we suppose that
$E_*$ is a parabolic bundle of rank two with the weights $\aa_i \leq 1 - \aa_i$ at $p_i$ and that $n \geq 1.$
Consider the following existence questions: 
\begin{enumerate}
\item[(I)] Does there exist $\Phi : E_* \rightarrow E_* \otimes K(D)$ with $(E_*, \Phi)$ stable? 
\item[(II)]  Does there exist $\Phi : E_* \rightarrow \widehat{E}_* \otimes K(D)$ with $(E_*, \Phi)$ stable?
\end{enumerate}

Such $\Phi$ are called stabilizing fields.
Of course, if $E_*$ is itself stable, then any $\Phi$ (e.g., $\Phi =0$)
gives us an affirmative answer. 
The other possibilities are if $E_*$ is unstable (meaning not
semistable) or if $E_*$ is strictly semistable. 
In either case, by choosing $L_*$ a line subbundle of maximal
parabolic degree, we get a short exact sequence
\begin{equation} \label{eq:ext1}
0 \longrightarrow L_* \stackrel{i}\longrightarrow E_* \stackrel{p}\longrightarrow M_* \longrightarrow 0
\end{equation}
with $\mu(L_*) \geq \mu(E_*).$
Let $\xi \in H^1(M^\vee_* \otimes L_*)$ be the extension class
representing (\ref{eq:ext1}).
If $E_*$ is unstable, then $\mu(L_*) > \mu (E_*)$ and
(\ref{eq:ext1}) is the Harder-Narasimhan filtration of $E_*$ and
is canonical.
If $E_*$ is
strictly semistable, then $\mu(L_*) = \mu(E_*)$ and (\ref{eq:ext1})
is the Jordan-H\"{o}lder filtration of $E_*$ and is not, in general,
canonical. For example, if $E_*$ is strictly semistable, then the subbundle 
$L_*$ is canonically determined if and only if the extension $\xi$ is nontrivial.

In the following proposition, the assumption $g \geq 2$ 
is not essential and after the proof, we treat the case $g \leq 1.$
\begin{prop}\label{prop:list-st}
If $g \geq 2$ and $E_*$ is 
not stable, then 
\begin{enumerate}
\item[(i)] $(E_*, \Phi)$ is a stable parabolic $K(D)$ pair for some $\Phi \Leftrightarrow
h^1(M^\vee_* \otimes \widehat{L}_*) \geq 1;$

\item[(ii)] $(E_*, \Phi)$ is a stable parabolic Higgs bundle for some $\Phi \Leftrightarrow
h^1(M^\vee_* \otimes L_*) > 1$ or $h^1(M^\vee_* \otimes L_*) = 1$ and $\xi = 0.$
\end{enumerate}
\end{prop}
\begin{pf}
Notice first of all that if such a $\Phi$ exists, then we can assume it is trace-free.
Now consider the short exact sequences of the sheaves of parabolic 
and strongly parabolic bundle endomorphisms
\begin{eqnarray} 
\label{eq:ses1}
&0 \rightarrow E_*^{\vee} \otimes L_* \otimes K(D)
\stackrel{\iota}{\longrightarrow} E_*^{\vee} \otimes_0 E_* \otimes K(D)
\stackrel{\pi}{\longrightarrow} L_*^{\vee} \otimes M_* \otimes K(D) \rightarrow 0,& \\
\label{eq:ses2}
&0 \rightarrow E_*^{\vee} \otimes \widehat{L}_* \otimes K(D)
\stackrel{\hat{\iota}}{\longrightarrow} E_*^{\vee} \otimes_0 \widehat{E}_* \otimes K(D)
\stackrel{\hat{\pi}}{\longrightarrow} L_*^{\vee} \otimes \widehat{M}_* \otimes K(D) \rightarrow 0,&
\end{eqnarray}
where $\pi, \hat{\pi}$ are the natural surjections, $\iota, \hat{\iota}$ are the natural isomorphisms to the kernels of $\pi, \hat{\pi}$ 
and $$E^\vee_* \otimes_0 E_* = \parend_0(E_*)$$ denotes the sheaf of
trace-free endomorphisms of $E_*.$ 
Notice that $H^0(E_*^{\vee} \otimes L_* \otimes K(D))$ 
and $H^0(E_*^{\vee} \otimes \widehat{L}_* \otimes K(D))$
are the relevant subspaces of fields
$\Phi$ for which $L_*$ is a $\Phi$-invariant subbundle.
If $(E_*, \Phi)$ is stable, then $L_*$ is not $\Phi$-invariant,
and $\pi_*(\Phi) \neq 0$ (similarly for $\hat{\pi}_*(\Phi)$).
This proves one implication of the following claim.
\begin{clm} \label{claim:1}
 Suppose that either $E_*$ is unstable or $\xi \neq 0,$ then
\begin{enumerate}
\item[(i)] for $\Phi \in H^0(E_*^{\vee} \otimes_0 E_* \otimes K(D)), (E_*, \Phi)$ is stable
$\Leftrightarrow 0 \neq \pi_*(\Phi);$
\item[(ii)]  for $\Phi \in H^0(E_*^{\vee} \otimes_0 \widehat{E}_* \otimes K(D)), (E_*, \Phi)$ is stable
$\Leftrightarrow 0 \neq \hat{\pi}_*(\Phi).$
\end{enumerate}
\end{clm}
\begin{pf}
To prove ($\Leftarrow$),
we just show that $L_*$ is the unique parabolic
subbundle of $E_*$ with $\mu(L_*) \geq \mu(E_*).$ Suppose
$L'_*$ is another such subbundle. If $E_*$ is unstable, then $\mu(E_*) > \mu(M_*)$ and
the projection $L'_* \rightarrow M_*$ is the zero map, which shows $L'_* = L_*$.
On the other hand, if $L'_* \rightarrow M_*$ is not the zero map, then it is an isomorphism
and defines a splitting of
(\ref{eq:ext1}), hence $\xi =0.$ 
\end{pf}

Now consider the coboundary maps in the cohomology sequences of (\ref{eq:ses1}) and (\ref{eq:ses2})
\begin{eqnarray*}
&H^0(L_*^{\vee} \otimes M_*\otimes K(D))
\stackrel{\delta}{\longrightarrow} H^1(E_*^{\vee}\otimes L_*\otimes K(D)),& \\
&H^0(L_*^{\vee} \otimes \widehat{M}_*\otimes K(D))
\stackrel{\hat{\delta}}{\longrightarrow} H^1(E_*^{\vee}\otimes \widehat{L}_*\otimes K(D)).&
\end{eqnarray*}
Here $\delta$ is the zero map since
by Serre duality
$$h^1(E_*^{\vee}\otimes L_*\otimes K(D)) = h^0 (L_*^{\vee} \otimes \widehat{E}_*) = h^0(L_*^{\vee} \otimes \widehat{L}_*) + h^0(L_*^{\vee} \otimes \widehat{M}_*)=0 .$$  
A diagram chase shows that the dual map of $\hat{\delta},$ $\hat{\delta}^\vee:H^0(L_*^\vee\otimes E_*)\rightarrow H^1(M_*^\vee\otimes L_*),$ maps $i$ to $\xi.$ Hence, 
$\hat{\delta}$ is the zero map if and only if $\xi =0.$ 
If $\xi \neq 0,$ then its image is one dimensional
because
$$h^1(E_*^{\vee}\otimes \widehat{L}_* \otimes K(D)) = h^0 (L_*^{\vee} \otimes E_*) = 
\begin{cases} 1 & \text{if $L_* \neq M_*$ or $\xi \neq 0,$} \\ 2 &\text{if $L_* = M_*$ and $\xi = 0.$} \end{cases}$$  
In the cases covered by the claim,
the proposition follows by another application of Serre duality 
$$h^0(L_*^{\vee}\otimes M_*\otimes K(D)) = h^1 (M_*^{\vee} \otimes \widehat{L}_*),$$  
$$h^0(L_*^{\vee}\otimes \widehat{M}_* \otimes K(D)) = h^1 (M_*^{\vee} \otimes L_*).$$  
The remaining cases follow by replacing the claim by the lemma below,
which we note is the only step of the argument where we use the assumption $g \geq 2.$
\end{pf}

\begin{lem} \label{lem:criterion}
If $g \geq 2$ and $E_*$ is not stable,
then
\begin{enumerate}
\item[(i)] $(E_*, \Phi)$ is a stable parabolic $K(D)$ pair for some
$\Phi \Leftrightarrow \ker \delta \neq 0;$
\item[(ii)] $(E_*, \Phi)$ is a stable parabolic Higgs bundle for some
$\Phi \Leftrightarrow \ker \hat{\delta} \neq 0.$
\end{enumerate}
\end{lem}
\begin{pf}
Since the lemma is a consequence of the claim, when it applies, we can 
assume that $E_*$ is strictly
semistable and $\xi =0.$
Furthermore, we only need to show ($\Leftarrow$).

We introduce some notation.
Define the intersection numbers $e_i$ and $\hat{e}_i$ by
\begin{eqnarray*} 
e_i &=& \begin{cases} \dim L_{p_i} \cap F_2(p_i) & \text{if $F_2(p_i) \neq 0$,} \\
				1 & \text{if $F_2(p_i) =0,$} \end{cases}\\
\hat{e}_i &=& \dim L_{p_i} \cap F_2(p_i).
\end{eqnarray*}
If 
$\beta_i = \hat{e}_i + (-1)^{\hat{e}_i} \alpha_i$ and $\gamma_i= 1-\beta_i$  
are the weights of $L_{p_i}$ and $M_{p_i},$ respectively, then
$$\hat{e}_i = \begin{cases}   0 & \text{if $\beta_i \leq \gamma_i$,}\\
                            1 & \text{if $\beta_i > \gamma_i$,}
\end{cases} \quad 
\hbox{ and }\quad
e_i =   
\begin{cases}   0 & \text{if $\beta_i < \gamma_i$,}\\
               1 & \text{if $\beta_i \geq \gamma_i$.}
\end{cases}$$
Set $|e| = \sum e_i$ and $|\hat{e}| = \sum \hat{e}_i$ and notice that 
$e_i > \beta_i - \gamma_i$ and  $\hat{e}_i \geq \beta_i - \gamma_i,$
with equality only when $\hat{e}_i = 0$ and $\beta_i = \gamma_i.$

If $\ker \delta \neq 0$ or $\ker \hat{\delta} \neq 0,$ then for generic $\Phi,$ $L_*$ is not $\Phi$-invariant.  
Suppose $L'_*$ ($\neq L_*$) is a line subbundle with $\mu(L'_*) \geq \mu(E_*).$ 
Semistability of $E_*$ implies $\mu(L'_*) = \mu(E_*).$
Then the restriction of $p$ to $L'_*,$ written $p_{L'} \, : \, L'_* \longrightarrow M_*,$
is an isomorphism since otherwise, $p_{L'} =0$
and $L'_* = L_*$.
Such subbundles are identified with sections of $p$ and are parameterized
by $H^0(M^\vee_* \otimes L_*).$ The relevant subspaces of $\Phi$ leaving
$L'_*$ invariant are $H^0(E^\vee_* \otimes M_* \otimes K(D))$
and $H^0(E^\vee_* \otimes \widehat{M}_* \otimes K(D)).$
Thus,
(i) will follow once we prove the inequality 
\begin{equation} \label{eq:ineq1}
h^0(M^\vee_* \otimes L_*) + h^0(E^\vee_* \otimes M_* \otimes K(D))< h^0(E^\vee_* \otimes_0 E_* \otimes K(D)),
\end{equation}
which is equivalent to $h^0(M^\vee_* \otimes L_*) < h^0(M^\vee_* \otimes L_* \otimes K(D)).$
Likewise, (ii) will follow from

\begin{equation} \label{eq:ineq2}
h^0(M^\vee_* \otimes L_*) + h^0(E^\vee_* \otimes \widehat{M}_* \otimes K(D))< h^0(E^\vee_* \otimes \widehat{E}_* \otimes K(D)),
\end{equation}
which is equivalent to $h^0(M^\vee_* \otimes L_*) < h^0(M^\vee_* \otimes \widehat{L}_* \otimes K(D)).$
Since $\mu(M^\vee_* \otimes L_*) = 0,$ 
$$h^0(M^\vee_* \otimes L_*) = 
\begin{cases} 0 & \text{if $M_* \neq L_*,$} \\
              1 & \text{if $M_* = L_*.$}
\end{cases}$$
On the other hand, because $h^1(M^\vee_* \otimes L_* \otimes K(D)) = h^0( L_*^\vee \otimes \widehat{M}_*) =0,$
it follows that 
\begin{eqnarray*}
h^0(M^\vee_* \otimes L_* \otimes K(D)) &=& \deg (M^\vee \otimes L \otimes K({\textstyle{\sum}^n_{i=1}} e_i p_i)) + \chi(X)\\
					&=& \deg L - \deg M +|e| + g-1.
\end{eqnarray*}
Notice that $\deg L - \deg M +|e| > \mu(L_*) - \mu(M_*) =0,$
hence (\ref{eq:ineq1}) holds provided
\begin{equation} \label{eq:cond1}
\deg L - \deg M +|e| \geq 2-g \hbox{ with equality } \Leftrightarrow L_* \neq M_*.
\end{equation}
This proves part (i) of the lemma when $g \geq 2.$
As for part (ii), notice that
$$h^1(M^\vee_* \otimes \widehat{L}_* \otimes K(D)) = h^0( L_*^\vee \otimes M_*) =
\begin{cases}  0 & \text{if $M_* \neq L_*,$} \\
              1 & \text{if $M_* = L_*,$}  
\end{cases}$$
and so (\ref{eq:ineq2}) follows as long as 
$\chi(M^\vee_* \otimes \widehat{L}_* \otimes K(D)) > 0.$  
We have
\begin{eqnarray*}
\chi (M^\vee_* \otimes \widehat{L}_* \otimes K(D)) &=& \deg(M^\vee \otimes L \otimes K({\textstyle{\sum}^n_{i=1}} \hat{e}_i p_i)) + \chi(X)\\
							&=& \deg L - \deg M +|\hat{e}| + g-1.
\end{eqnarray*} 
Hence (\ref{eq:ineq2}) holds provided
\begin{equation} \label{eq:cond2}
\deg L - \deg M +|\hat{e}| \geq 2-g.
\end{equation}
But $\deg L - \deg M +|\hat{e}| \geq \mu(L_*) - \mu(M_*) = 0$
(with equality implying that $\beta_i = \gamma_i$ for all $i$).
This proves part (ii) of the lemma when $g \geq 2.$
\end{pf}


One can deduce the following corollary using Riemann-Roch.
\begin{cor}
If $g \geq 3,$ then for every semistable $E_*$, there exists
a Higgs field $\Phi$ making $(E_*, \Phi)$ a stable parabolic Higgs bundle. 
\end{cor}

We now explain how to extend these results to lower genus. 
Clearly, the proposition holds for $g \leq 1$ whenever $E_*$ is unstable or $\xi \neq 0$
by virtue of the claim. 
So assume that $E_*$ is semistable and $\xi = 0.$
The only place 
where we make essential use of the assumption $g \geq 2$ is 
in the proof of Lemma \ref{lem:criterion}.
In particular, we observe from (\ref{eq:cond1}) and (\ref{eq:cond2}) that the 
inequalities 
(\ref{eq:ineq1}) and (\ref{eq:ineq2})
fail (respectively) if
\begin{enumerate}
\item[(i)] $0 < \deg L - \deg M +|e| \leq 2-g$  with equality $ \Leftrightarrow L_* = M_*,$
\item[(ii)] $0 \leq \deg L - \deg M +|\hat{e}| \leq 1-g.$  
\end{enumerate}

Thus, the only counterexamples to Lemma \ref{lem:criterion} for $g \leq 1$
are given by the semistable, split bundles $E_*$ satisfying (i) and (ii) 
along with the additional requirements
(${\hbox{\rm i}}'$) $\ker \delta \neq 0$ and (${\hbox{\rm ii}}'$) $\ker \hat{\delta} \neq 0.$ 
First, we list these counterexamples to Lemma \ref{lem:criterion},
then we show that the bundles satisfying (i) and (ii) never give
rise to any stable parabolic $K(D)$ pairs or stable parabolic Higgs bundles, respectively.

If $E_*$ is semistable and split and satisfies 
(i) and (${\hbox{\rm i}}'$), i.e., if $h^1(M^\vee_* \otimes \widehat{L}_*) \geq 1,$
then there are but two possibilities:
\begin{enumerate}
\item[(i--a)] $(g,n) \in \{(0,2), (1,1)\},\, E_*= L_* \oplus M_*$ and $L_* = M_*,$
\item[(i--b)] $g=0,\, E_*= L_* \oplus M_*, \, \mu(L_*) = \mu(M_*), \, \deg L - \deg M + |e| = 1.$ 
\end{enumerate}
Now if $E_*$ is semistable and split and satisfies (ii) and (${\hbox{\rm ii}}''$), 
i.e., if $h^1(M^\vee_* \otimes L_*) \geq 1,$ 
then again, we have only two possibilities:
\begin{enumerate}
\item[(ii--a)] $g=0,\,  E_*= L_* \oplus M_*,\, \mu(L_*) = \mu(M_*),$ and $ 0 \leq \deg L - \deg M +|\hat{e}| \leq  1,$
\item[(ii--b)] $g=1,\, E_*= L_* \oplus M_*,\, L_* = M_*.$ 
\end{enumerate} 

We now show that if $E_*$ satisfies (i), then $(E_*, \Phi)$ is not stable
for any $\Phi \in H^0(E^\vee_* \otimes_0 E_* \otimes K(D))$ and
if $E_*$ satisfies (ii),
then $(E_*, \Phi)$ is not stable
for any $\Phi \in H^0(E^\vee_* \otimes_0 \widehat{E}_* \otimes K(D)).$ 
For example, suppose that $\deg L - \deg M +|e| = 2-g$ in (i), so that $L_* = M_*.$
Then either $g=0$ and $n=2$ or $g=1=n.$ In either case,
$$H^0(E^\vee_* \otimes_0 E_* \otimes K(D)) = H^0(K(D))^{\oplus 3} = H^0({\cal O}_X)^{\oplus 3}.$$
Thus, any $\Phi$ 
is a constant matrix, one of whose eigenspaces determines a $\Phi$-invariant subbundle
violating the condition for stability.
Otherwise, if $\deg L - \deg M +|e| = 1-g$ in (i), then $h^0(M^\vee \otimes L_* \otimes K(D)) = 0$
so that $M_*$ is $\Phi$-invariant for all $\Phi.$

As for (ii), suppose first of all that $g=0$ and $\deg L - \deg M +|\hat{e}| \leq 1.$
Then $h^0(M^\vee_* \otimes \widehat{L}_* \otimes K(D)) =0$ and $M_*$ is $\Phi$-invariant for all $\Phi.$
Now if $g=1$ and $\deg L - \deg M +|\hat{e}| =0,$ 
then either $L_* \neq M_*$ and $M_*$ is $\Phi$-invariant for all $\Phi$ 
or $L_* = M_*$ and $H^0(E^\vee_* \otimes_0 \widehat{E}_* \otimes K(D)) = H^0({\cal O}_X)^{\oplus 3},$
in which case every $\Phi$ is a constant matrix.
This proves the following proposition.
\begin{prop} \label{prop:crit}
If $E_*$ is not stable and $g \leq 1,$ then
\begin{enumerate}
\item[(i)] $(E_*, \Phi)$ is a stable parabolic $K(D)$ pair for some $\Phi \Leftrightarrow \,
E_*$ is not one of the bundles occurring in \rom(i--a\rom) or \rom(i--b\rom) and $h^1(M^\vee_* \otimes \widehat{L}_*) \geq 1;$ 
\item[(ii)] $(E_*, \Phi)$ is a stable parabolic Higgs bundle for some $\Phi \Leftrightarrow
\, E_*$ is not one of the bundles occurring in \rom(ii--a\rom) or \rom(ii--b\rom) and 
either $h^1(M^\vee_* \otimes L_*) > 1$ or $h^1(M^\vee_* \otimes L_*) = 1$ and $\xi = 0.$
\end{enumerate}
\end{prop}

We could ask questions (I) and (II) replacing stability with semistability.
Of course, if $E_*$ itself is semistable, then so is $(E_*, \Phi)$ for any $\Phi.$
So we can assume that $E_*$ is unstable and apply the claim to determine 
precisely which $\Phi$ make $(E_*, \Phi)$ stable. 
One last comment is that if $(E_*, \Phi)$ is strictly semistable, then $E_*$ 
must also be strictly semistable. The converse, however, is false.

\subsection{Example: Rank 2 parabolic bundles over ${\Bbb P}^1$ with 3 parabolic points}	

In this section, we describe the moduli spaces
${\cal M}_\alpha, {\cal N}_\alpha$ and ${\cal P}_\alpha$ of rank two
bundles over $X={\Bbb P}^1$
with parabolic points in the reduced divisor $D=p_1+p_2+p_3.$

This case seems trivial as it turns out that
${\cal N}_\alpha$ is always just one point
and that ${\cal P}_\alpha$ is always just the affine space ${\Bbb C}^5$.
However, our complete description of this case
sheds light on the general phenomenon that the moduli
spaces ${\cal N}_\alpha$ and ${\cal P}_\alpha$ do not change when
the weights are permitted to vary
(even when ${\cal M}_\alpha$ becomes empty!).
This trivial case is a prototype for such behavior.
    
The simplest nontrivial cases are
$X={\Bbb P}^1$ with 4 parabolic points
and $X=C,$ an elliptic curve, with one parabolic point.
In either case, ${\cal M}_\alpha,$ if nonempty, is ${\Bbb P}^1,$
and ${\cal N}_\alpha$ is a connected
nonsingular noncompact surface containing the cotangent bundle of
${\Bbb P}^1$. There is a proper map from ${\cal N}_\alpha$ to 
${\Bbb C}$ called the Hitchin map whose fibers over nonzero points 
$t \in {\Bbb C}$ are elliptic curves and whose fiber over $0$ is 
a union of five
rational curves arranged in a $\widetilde{D}_4$ configuration. 
This case will be treated in the second part of this paper.
	       
We suppose that $\mu(E_*)=0$ and that the weights at $p_i$ are $\alpha_i$ and $1-\alpha_i$ 
for some $\alpha \in W = \{ (\alpha_1,\alpha_2,\alpha_3) \mid 0 < \alpha_i < \frac{1}{2} \}.$
Note that this is equivalent to saying that $\det E_* = {\cal O}_X$ (as parabolic bundles) and $E_*$ has full flags at each $p_i.$
For $e=(e_1,e_2,e_3),$ where $e_i \in \{0,1\},$ we use $\beta(\alpha,e)$ (or simply $\beta$)
to denote the weights $\beta_i = e_i + (-1)^{e_i} \alpha_i.$
Let $$I = \{ (0,0,0), (0,1,1),(1,0,1),(1,1,0) \}.$$
Inside $W$ there are four hyperplanes 
$$H_e = \{ \alpha \mid \beta(\alpha,e) = 1 + \frac{|e|}{2} \}$$ for
$e \in I$ whose complement $W \setminus \bigcup_{e \in I} H_e$ consists of five chambers:
$C_e = \{ \alpha  \mid \beta(\alpha,e) > 1 + \frac{|e|}{2} \}$
for $e \in I$ and $C_0  = \{ \alpha  \mid \beta(\alpha,e) < 1 + \frac{|e|}{2} \hbox{ for all }  e\in I\}.$

The following is an immediate consequence of the criteria established in the previous section.
\begin{lem}
If $(E_*, \Phi)$ is a semistable $K(D)$ pair, then
the bundle $E_*$ is described as an extension 
\begin{equation} \label{eqn:ext2}
0 \longrightarrow L_* \longrightarrow E_* \longrightarrow L_*^\vee \longrightarrow 0
\end{equation}
where $L_*$ satisfies  $h^1(L_*^{\otimes 2})=1.$
\end{lem}
\begin{pf}
If $E_*$ is not stable, then by Proposition \ref{prop:crit}, we see that
$h^1(L_*^{\otimes 2}) \geq 1.$ Since $\mu(L_*) \geq 0,$ we see that  $\mu(L_*^{\otimes 2}) \geq 0,$
and because there are only three weights, this implies
$\deg(L_*^{\otimes 2})_0 \geq -2.$ Thus $h^1(L_*^{\otimes 2}) = 1.$

If $E_*$ is stable, then
by Grothendieck's Theorem, $E = {\cal O}(-1) \oplus {\cal O}(-2).$
Let $L_*$ be ${\cal O}(-1)$ with weights inherited as a subbundle of $E_*.$
Notice that $h^1(L_*^{\otimes 2}) \leq 1.$
But by stability of $E_*,$ the extension (\ref{eqn:ext2}) must be nontrivial,
so $h^1(L_*^{\otimes 2}) \geq 1.$ 
\end{pf}

We now determine all possible line subbundles $L_*$ with $h^1(L_*^{\otimes 2})=1.$
For fixed $\aa \in W,$ there are four 
possible line subbundles $L_*$ with $h^1(L_*^{\otimes 2})=1,$ namely
$${L^e_*} = {\cal O}_X(-1-\frac{|e|}{2})[-\beta(\alpha,e)]$$ for $e \in I.$ 
We denote by $G^e_*$ the nontrivial extension gotten from (\ref{eqn:ext2}) with $L_*=L^e_*$.
Notice that $G^e_*$ is unique up to isomorphism because $h^1({L^e_*}^{\otimes 2})=1.$ 
Let $F^e_*=L^e_*\oplus {L^e_*}^\vee.$

It is not hard to see that $G^e_*$ and $G^{e'}_*$ are isomorphic
for $e, e' \in I.$ Set $G_* = G^e_*.$
This, together with the previous lemma, shows that if $(E_*, \Phi)$ is semistable,
then $E_*$ is one of the five bundles in the set $\{ G_*, F^e_* \}.$
 
Recall that two bundles $E_*$ and $E'_*$ are called S-equivalent (written $E_* \sim_S E'_*$)
if their associated graded bundles are isomorphic, 
i.e., if $\gr E_* \simeq \gr E'_*.$
We use $E_*$ to denote the isomorphism class of a bundle and $[E_*]$ for its S-equivalence class.
\begin{prop} 
\begin{enumerate}
\item If $\aa \in C_0,$ then  ${\cal M}_\aa = \{{G_*}\}.$ 
\item If $\aa \in C_e, e \in I,$ then ${\cal M}_\aa = \emptyset.$ 
\item If $\aa \in H_e,$ then ${\cal M}_\aa = \{ [F^e_*] \}$ and ${G_*} \sim_S F^e_*$ are the two 
distinct isomorphism classes of semistable bundles.
\end{enumerate}
\end{prop}
\begin{pf}
From the above considerations, if $E_*$ is semistable, then $E_* =G_*$ or $F^e_*.$
But $G_*$ is stable if and only if $\aa \in C_0,$ and $F^e_*$ is never stable.
On the other hand, 
if $\aa \in H_e,$ then $G_*$ and $F^e_*$ are clearly strictly semistable 
with associated graded bundle $L^e_* \oplus (L^e_*)^\vee.$ 
\end{pf}

The next lemma shows which auxiliary fields can arise for these five bundles.
\begin{lem}
For any $\aa \in W,$ we have
\begin{enumerate}  
\item[(i)] $G_*$ is simple,  
$h^0(G^\vee_* \otimes G_* \otimes K(D)) = 5,$ and $h^0(G^\vee_* \otimes \widehat{G}_* \otimes K(D))=0,$
\item[(ii)] $\Aut F^e_* = {\Bbb C}^*\! \times {\Bbb C}^*,  \, h^0({F^e_*}^\vee \otimes F^e_* \otimes K(D)) = 5,$
and $h^0({F^e_*}^\vee \otimes \widehat{F}^e_* \otimes K(D))=1.$ 
\end{enumerate}
\end{lem}
\begin{pf}
For $\aa \in C_0, \, G_*$ is stable, and therefore simple. But this 
property is independent of the weights,
and it follows that for any $\aa \in W,$ 
\begin{eqnarray*}
1 &=& h^0(G_*^\vee\otimes G_*) = h^1(G_*^\vee\otimes \widehat{G}_* \otimes K(D)), \\
0 &=& h^0(G_*^\vee\otimes \widehat{G}_*) = h^1(G_*^\vee\otimes{G}_* \otimes K(D)).
\end{eqnarray*}
Direct computation shows $\deg (G^\vee_* \otimes G_*)_0 = -3$
and $\deg (G^\vee_* \otimes \widehat{G}_*)_0 = -9,$
and part (i) follows using  $K(D) = {\cal O}(1)$ and Riemann-Roch.

As for (ii), since $({L^e_*}^{\otimes 2})_0 = {\cal O}(-2)$  and
$({L^e_*}^{\otimes -2})_0 = {\cal O}(-1),$
every automorphism of $F^e_*$ is diagonal and $\Aut {F^e_*}= {\Bbb C}^*\! \times {\Bbb C}^*.$ 
Also, $h^0({L^e_*}^{\otimes 2}\otimes K(D))= 0,$
so every $\Phi\in H^0({F^e_*}^\vee\otimes F^e_*\otimes K(D))$ has the form
$$ \Phi= \begin{pmatrix} a_1 & 0 \\ \phi & a_2 \end{pmatrix} $$
with $\phi\in H^0(L_*^{\otimes -2}\otimes K(D))={\Bbb C}$ and $a_i \in H^0(K(D))={\Bbb C}^2.$
Morevoer, $H^0({F^e_*}^\vee\otimes \widehat{F}^e_*\otimes K(D)) = H^0(L_*^{\otimes -2}\otimes K(D)),$
which completes the proof of part (ii)
\end{pf}

We can identify
the action of $\Aut(F^e_*)$ on $H^0({F^e_*}^\vee\otimes F^e_*\otimes K(D)),$ 
it is given by conjugation
$$(z_1,z_2) \cdot \Phi = 
\begin{pmatrix} z_1 & 0 \\0& z_2 \end{pmatrix}
 \begin{pmatrix} a_1 & 0 \\ \phi & a_2 \end{pmatrix}
 \begin{pmatrix} z_1^{-1} & 0 \\0& z_2^{-1} \end{pmatrix} =
\begin{pmatrix} a_1 & 0 \\z_1^{-1}z_2 \phi & a_2 \end{pmatrix}.$$

Suppose that $\aa\in C_e$ and set $V=\Ext^1({L^e_*}^\vee,L^e_*)={\Bbb C}$.
Let ${\cal E}_*$ be the universal parabolic bundle on $V \times X$ which, when restricted to
$\{\xi\} \times X,$ is the bundle $G^\xi_*$ in (\ref{eqn:ext2}) with 
$L_* = L^e_*$ and extension class $\xi.$ For $\xi \neq 0, \, G^\xi_* \simeq G_*$ and 
obviously $G^0_* = F^e_*.$

Let $p_X$ and $p_V$ denote the two projection maps from $V \times X$
and define ${\cal L}_*$ to be the pullback bundle $p_X^* L^e_*.$
Consider the direct image sheaves of ${\cal E}^\vee_* \otimes {\cal E}_* \otimes K(D)$
and ${\cal L}_*^{\otimes -2} \otimes K(D)$ under $p_V,$ which, by the previous lemma,
are locally free sheaves over
$V$ whose associated vector bundles, $M$ and $N,$ are trivial with ranks  5 and 1, respectively.
Notice that
$N$ is canonically isomorphic to $V \times H^0({L^e_*}^{\otimes -2} \otimes K(D)).$
This is key to following construction.

The canonical map 
$\tilde{\pi} : {\cal E}^\vee_* \otimes {\cal E}_* \otimes K(D) \rightarrow {\cal L}_*^{\otimes -2} \otimes K(D)$
of the previous section induces 
$ \tilde{\pi}_* : M \rightarrow N$ which is surjective,
because the restriction of $\tilde{\pi}_*$
to a fiber above $\xi$ can be identified with 
$\pi^\xi_* : H^0({G^\xi_*}^\vee\otimes G^\xi_* \otimes K(D))\rightarrow H^0({L^e_*}^{\otimes -2}\otimes K(D)),$
whose cokernel is $H^1({G^\xi_*}^\vee \otimes L^e_* \otimes K(D)) = 0.$
Fix some $0 \neq \phi_0 \in H^0({L^e_*}^{\otimes -2}\otimes K(D))$ and set
$Y =  \tilde{\pi}_*^{-1}( \{\phi_0\} \times V) \simeq {\Bbb C}^4 \times V.$
\begin{prop}
\begin{enumerate}
\item If $\aa\in C_0,$ then ${\cal P}_\aa \simeq H^0(G_*^\vee\otimes G_* \otimes K(D)) = {\Bbb C}^5$
and ${\cal P}_\aa^0 \simeq {\Bbb C}^3$. 
\item If $\aa\in C_e,$ then ${\cal P}_\aa\simeq Y\simeq {\Bbb C}^5$ and ${\cal P}_\aa^0 = {\Bbb C}^3.$ 
\item If $\aa\in H_e,$ then ${\cal P}_\aa \simeq H^0(G_*^\vee\otimes {G_*} \otimes K(D)) = {\Bbb C}^5$
and its strictly semistable part can be identified with a hyperplane.
\end{enumerate}
\end{prop}
\noindent
{\it Remark.\,} In the course of the proof, we will
determine the isomorphism classes of semistable parabolic $K(D)$ pairs.
This differs from the above only for strictly semistable bundles,
because the S-equivalence class of a stable bundle
is precisely its isomorphism class.
For $\aa \in H_e,$ we will find that there are three distinct components of
isomorphism classes of strictly semistable bundles, each is just a copy of
${\Bbb C}^4.$ 
\begin{pf}
Part (1) follows from the fact that $h^0({L^e_*}^{\otimes 2}\otimes K(D))=0,$ hence 
${L^e_*}^\vee$ is a $\Phi$-invariant subbundle of $F^e_*$ for any $\Phi.$ 
Thus, if $(E_*, \Phi)$ is stable and $\aa \in C_0,$ then $E_* \simeq G_*.$
For part (2), 
if $(\xi,\Phi) \in Y,$ then the associated $K(D)$ pair $(G^\xi_*, \Phi)$ 
is stable by Claim \ref{claim:1}
since $\pi^\xi_* (\Phi) = \phi_0 \neq 0$.
This gives a map $\eta:Y\rightarrow {\cal P}_\aa,$
which we claim is a bijection. 
To see this, write $Y = Y' \cup Y'',$
where $Y' = Y|_{V \setminus 0}$ and $Y'' = Y|_0,$
and ${\cal P}_\aa = {\cal P}_\aa' \cup {\cal P}_\aa'',$
where ${\cal P}_\aa'$ and ${\cal P}_\aa''$ consist of the $K(D)$ pairs $(E_*, \Phi)$
with underlying bundle $E_*$ isomorphic to  $G_*$ and $F^e_*,$ respectively.
The restriction of $M$ to $V\setminus\{0\}$ is naturally isomorphic to 
$(V\setminus\{0\})\times H^0(G_*^\vee\otimes G_*\otimes K(D)).$
For $(\xi,\Phi)\in M|_{V\setminus\{0\}}$ and $t\in{\Bbb C}^*,$
        $$\pi^{t \xi}_*(\Phi)=t^{-1}\pi^\xi_*(\Phi).$$
It follows from this formula that $\eta$ induces a bijection between $Y'$ and ${\cal P}_\aa'$

Using the description of the action of $\Aut F^e_*$ on $H^0({F^e_*}^\vee \otimes F^e_* \otimes K(D))$
following the proof of the previous lemma,
every $(F^e_*, \Phi) \in {\cal P}_\aa''$ is isomorphic to
$(F^e_*, \Phi_0)$ where
$$ \Phi_0 = \begin{pmatrix} a_1 & 0 \\ \phi_0 & a_2 \end{pmatrix}. $$
Hence, $\eta$ gives a bijection between $Y''$ and ${\cal P}_\aa''.$ 

To prove (3), notice that we have a map
$\eta:H^0(G_*^\vee\otimes {G_*}\otimes K(D))\rightarrow {\cal P}_\aa.$ 
Since ${\cal P}_\aa$ is normal, it is enough to show that this map is bijective. 
Now by Claim \ref{claim:1}, we see that 
${\cal P}_\aa^s\simeq H^0(G_*^\vee\otimes {G_*}\otimes K(D))\setminus\Ker\pi_*.$ 
The strictly semistable bundles are pairs of the form $(F^e_*,\Phi)$ for
any $\Phi,$ and
$({G_*},\Phi)$ with $\Phi\in\Ker\pi_*.$ 
If $\Phi\in\Ker\pi_*$, 
the subbundle $L^e_*$ is $\Phi$-invariant and we get the extension
of parabolic $K(D)$ pairs
\begin{equation}\label{eq:ext}
	0 \longrightarrow (L^e_*,\phi)
	\longrightarrow (G_*,\Phi)
	\longrightarrow ({L^e_*}^\vee,\psi)
	\longrightarrow 0.
\end{equation}
Thus $\gr (G_*, \Phi) = (L^e_*,\phi) \oplus ({L^e_*}^\vee,\psi)$ for $\Phi \in\Ker\pi_*.$
Consider now the map
$$\lambda:\Ker\pi_* \longrightarrow H^0({L^e_*}^\vee\otimes 
L^e_*\otimes K(D))\oplus H^0(L^e_*\otimes {L^e_*}^\vee\otimes K(D))$$
defined by $\Phi\mapsto(\phi,\psi)$.
For $\phi=\psi=0$, then the extension 
(\ref{eq:ext}) induces the zero map ${L^e_*}^\vee \rightarrow L^e_*\otimes K(D)$ 
(because $H^0({L^e_*}^{\otimes 2} \otimes K(D))=0$) 
and it follows that $\Phi=0.$  
So $\lambda$ is injective. But the domain and range of
$\lambda$ are both 4-dimensional, and so $\lambda$ is an isomorphism.
Clearly $\gr(F^e_*,\Phi) = (L^e_*,\phi)\oplus({L^e_*}^\vee,\psi),$ and it follows that
$\lambda$ gives a bijection between $\Ker\pi_*$ and ${\cal P}_\aa^{sss}.$
\end{pf}

Choosing some $0 \neq \Phi_0 \in H^0({F^e_*}^\vee \otimes \widehat{F^e_*} \otimes K(D)) = {\Bbb C}$ and
using the action of $\Aut(F^e_*),$ 
it is easy to verify that $(F^e_*, \Phi)$ is isomorphic to $(F^e_*,  \Phi_0)$ for all $\Phi \neq 0.$
The proof of the last proposition is left as an entertaining exercise in applying the above lemmas.

\begin{prop}
\begin{enumerate}
\item If $\aa \in C_0,$ then  ${\cal N}_\aa = \{({G_*},0)\}.$ 
\item If $\aa \in C_e, e \in I,$ then ${\cal N}_\aa = \{(F^e_*,\Phi_0)\}.$ 
\item If $\aa \in H_e,$ then ${\cal N}_\aa = \{ [F^e_*,0] \}$ and $(G_*,0) \sim_S (F^e_*,0) \sim_S (F^e_*,\Phi_0)$
are the three distinct isomorphism classes of semistable Higgs bundles.
\end{enumerate}
\end{prop}

\section{A Topological Description of ${\cal N}^0_\alpha$ in Rank Two}
\subsection{The function spaces of Biquard and construction of Konno}

We begin with a brief overview of the gauge theoretical description of
${\cal N}_\alpha$ following \cite{konno}.

It is convenient to think of the parabolic bundle separate from its holomorphic structure,
so we use $E_*$ to denote the underlying topological parabolic bundle (weights $\alpha$)
and $\delbar_E$ its holomorphic structure.
By tensoring with an appropriate line bundle,
we can always assume that $\mu(E_*) = 0.$ 
We shall also restrict our attention to generic weights,
i.e., weights $\aa$ for which $\alpha$-stability and $\alpha$-semistability coincide.
Let ${\cal C}$ denote the affine space of all holomorphic structures on $E,$
and ${\cal G}_{\Bbb C}$ the group of smooth bundle automorphisms of $E$
preserving the flag structure.
Introduce a metric $\kappa$ adapted to $E$ 
($\kappa$ is unitary and smooth on $E|_{X \setminus D}$, but singular at $p \in D$ in a prescribed way,
see Definition 2.3 \cite{biquard}),
and let ${\cal A}$ denote the affine space of $\kappa$-unitary connections. 
Define ${\cal G}$ to be the subgroup of ${\cal G}_{\Bbb C}$
consisting of $\kappa$-unitary gauge transformations.
Letting ${\cal C}_{ss}$  and ${\cal A}_{\hbox{\scriptsize \it flat}}$ be
the subspaces of
$\alpha$-semistable holomorphic structures and the flat connections, respectively,
Biquard proved that 
$${\cal M}_\alpha \stackrel{\hbox{\scriptsize def}}{=} {\cal C}_{ss} / {\cal G}_{\Bbb C} \cong 
{\cal A}_{\hbox{\scriptsize \it flat}}/{\cal G}$$
by introducing the norms 
$\| \; \|_{D^p_k},$ defining the weighted Sobolev spaces
${\cal C}^p$ and ${\cal A}^p$ of $D^p_1$
holomorphic structures and $D^p_1 \, \kappa$-unitary connections, and taking quotients by
the groups
${\cal G}_{\Bbb C}^p$ and ${\cal G}^p$  of $D^p_2$ gauge transformations
for a certain $p>1$ \cite{biquard}.

The same approach works for parabolic Higgs moduli, at least for generic weights,
as was shown by Konno. The arguments in \cite{konno} are given for
moduli with fixed determinant, but remain equally valid without this condition.
We set
\begin{eqnarray*}
{\cal H} &=& \{ (\delbar_E, \Phi) \in {\cal C} \times \Omega^{1,0}(\End E) \mid
\delbar_E \Phi = 0 \hbox{ on } X \setminus D \hbox{ and at each } p \in D, \\
& & \Phi \hbox{ has a simple pole with nilpotent residue with respect to the flag} \}.
\end{eqnarray*}
Note that ${\cal H}$ (this is denoted by ${\cal D}$ in \cite{konno})
is just the differential geometric definition of the space of parabolic Higgs bundle structures
on $E_*,$ for example, the nilpotency condition implies that $\Phi$ is strongly parabolic.

For $A \in {\cal A},$ we use $d_A$ for its covariant derivative, $F_A$ for its curvature, 
and $d_A''$ for the $(0,1)$ component of $d_A,$ so $d_A'' \in {\cal C}.$
Define ${\cal E} = {\cal A} \times \Omega^{0,1}(\End E)$ and ${\cal E}^p$
as its completion with respect to the norms $\| \; \|_{D^p_1},$ 
and set
$${\cal E}_{\hbox{\scriptsize \it flat}} = \{ (d_A, \Phi) \in {\cal E}^p \mid d_A'' \Phi = 0, F_A + [\Phi, \Phi^*] = 0 \}.$$ 
(This last space is denoted ${\cal D}^p_{H \! E}$ by Konno.)
Using the usual definition of stability on ${\cal H},$ 
Theorem 1.6 of \cite{konno} shows that for some $p>1,$
$${\cal N}_\alpha \stackrel{\hbox{\scriptsize def}}{=} {\cal H}_{ss} / {\cal G}_{\Bbb C} 
\cong {\cal E}_{\hbox{\scriptsize \it flat}}/{\cal G}^p.$$
The advantage of the second quotient 
is that it endows ${\cal N}_\alpha$ with a natural hyperk\"ahler structure, namely
by viewing it as a hyperk\"ahler quotient of ${\cal E}^p$ (in the sense of \cite{HKLR}), 
whose hyperk\"ahler structure is
given by the metric
$$g((\xi,\phi),(\xi,\phi)) = 2 i \int_X \Tr(\xi^* \xi + \phi \phi^*),$$ 
which is K\"ahler with respect to each of 
three complex structures
$$I(\xi,\phi) = (i \xi, i \phi), \quad
J(\xi,\phi) = (i \phi^*,- i \xi^*), \quad
K(\xi,\phi) = (- \phi^*, \xi^*).$$

\subsection{The Morse function for the moduli space of parabolic Higgs bundles}
Assume that $E_*$ is a rank two parabolic bundle  
with generic weights $\alpha_i$ and $1-\alpha_i$ at $p_i$ and that $\mu_\alpha(E_*)=0.$
Write $\alpha=(\alpha_1,\ldots,\alpha_n).$ We will always assume $n\ge 1.$
We consider the moduli with fixed determinant and trace-free Higgs fields,
requiring the following minor modifications in the definitions of the previous section:
\begin{itemize}
\item[(i)] the induced connection $d_\Lambda$ or holomorphic structure
$\delbar_\Lambda$ on $\Lambda^2 E$ be fixed; 
\item[(ii)] the Higgs field be trace-free, i.e. $\Phi \in \Omega^{1,0}(\End_0 E).$
\end{itemize}
We denote the corresponding spaces by ${\cal A}^0, {\cal C}^0, {\cal E}^0,$ and ${\cal H}^0.$

As in \cite{hitchin}, we consider the circle action defined on ${\cal E}^0$ by
$e^{i\theta} \cdot (d_A,\Phi) = (d_A, e^{i\theta}\Phi).$
This action preserves the subspace ${\cal E}^0_{\hbox{\scriptsize \it flat}}$ and commutes with
the action of the gauge group ${\cal G}^p,$ thus it
descends to give a circle action $\rho$ on ${\cal N}^0_\alpha.$
This action 
commutes with the complex structure defined by $I$
and preserves the symplectic form $\omega_1(X,Y)=g(IX,Y),$
so the associated moment map $\mu_\rho(d_A, \Phi) = \frac{1}{4 \pi}\| \Phi \|_{D^p_1}^2,$
renormalized for convenience, is a Bott-Morse function and can be used
to determine the Betti numbers of ${\cal N}^0_\alpha.$

We introduce some notation which will be used throughout the rest of this section.
For any line subbundle $L_*$ of $E_*,$
let $e_i(L) = \dim L_{p_i} \cap F_2(p_i) \in \{0,1\}.$ The weight inherited by $L_*$ is then
$\beta_i(L) = e_i +(-1)^{e_i} \alpha_i.$ We will often suppress the dependence on
$L$ and simply write $e = (e_1, \ldots, e_n)$ and $\beta=(\beta_1, \ldots, \beta_n).$
We will also write $\beta(\alpha,e)$ when 
we want to emphasize the functional
dependence of $\beta$ on $\alpha$ and $e.$ 
We also use $|e| = \sum_{i=1}^n e_i.$

\begin{thm} \label{thm:morse}
\begin{itemize}
\item[(a)] The map $\mu_\rho :{\cal N}^0_\alpha \longrightarrow {\Bbb R}$ is 
a proper Morse function.
\item[(b)] Whenever nonempty, ${\cal M}^0_\alpha$ is the unique critical submanifold corresponding to
the minimum value $\mu_\rho =0.$
The other critical submanifolds are given by ${\cal M}_{d,e}$ for an integer $d$ 
and $e \in {\Bbb Z}_2^n$ satisfying
\begin{equation} \label{3a}
-\sum_{i=1}^n \beta_i(\alpha,e) < d \leq g-1 -|e|/2.
\end{equation}
Along ${\cal M}_{d,e}, \, \mu_\rho$ takes the value  $d + \sum_{i=1}^n \beta_i.$
\item[(c)] The critical submanifold ${\cal M}_{d,e}$ is 
$\widetilde{S}^{h_{d,e}}X,$ 
the $2^{2g}$ cover of the symmetric product 
$S^{h_{d,e}}X$  under the map $x \mapsto 2x$ on $J_X.$ Here, $h_{d,e} = 2g -2 -2d - |e|.$ 
\item[(d)] The Morse index of ${\cal M}_{d,e}$ is given by 
$\lambda_{d,e} = 2(n+2d +g-1 +|e|).$
\end{itemize}
\end{thm}
\noindent
{\it Remark. \,} If $g=0,$ there are always $\alpha$ with ${\cal M}^0_\alpha = \emptyset$
(but ${\cal N}^0_\alpha \neq \emptyset$).
For these $\alpha,$ the minimum value is achieved along some ${\cal M}_{d,e},$
which we identify in the next section.

\begin{pf}
Properness of $\mu_\rho$ follows from the global compactness result for parabolic bundles
of Biquard (Theorem 2.14 in \cite{biquard}). This proves (a).
All the other statements rely on 
the following correspondence between the circle action and
the moment map given in \cite{frankel}. 
\begin{itemize}
\item[(1)] Critical submanifolds are connected components of the fixed point set of $\rho.$ 
\item[(2)] The Morse index of a critical submanifold
equals the dimension of the negative weight
space of the infinitesimal circle action on its normal bundle.
\end{itemize}
Suppose that $(d_A,\Phi)$ is a fixed point of the circle action upstairs in ${\cal E}_{\hbox{\scriptsize \it flat}}.$
Then $\Phi = 0$ and this shows that one component of the fixed point set in ${\cal N}^0_\alpha$
consists of ${\cal M}^0_\alpha,$ the moduli of stable parabolic bundles with fixed determinant.

The other fixed points arise from when
$e^{i \theta}\cdot (d_A, \Phi)$ is gauge equivalent to $(d_A,\Phi),$
i.e., when there is a one parameter family $g_\theta \in {\cal G}^p$ such that
\begin{eqnarray*}
g^{-1}_\theta \Phi g_\theta &=& e^{i\theta}\Phi,\\
g^{-1}_\theta d_A  g_\theta &=& d_A.
\end{eqnarray*}
By the first equation, $g_\theta$ is not central, and
by the second, we see that $d_A$ is reducible and consequently
the holomorphic parabolic bundle splits according to the
eigenvalues of $g_\theta.$
Write $E_* = L_* \oplus M_*$ as a direct sum of parabolic bundles.
We assume (wlog) that $\mu_\alpha(L_*) > 0 > \mu_\alpha(M_*).$
Let $d=\deg L$ and $e=(e_1,\ldots,e_n)$ where $e_i = \dim L_{p_i} \cap F_2(p_i).$ 
Then $L$ inherits the weight $\beta_i = e_i + (-1)^{e_i} \alpha_i$ at $p_i$
as a parabolic subbundle of $E_*$
and 
\begin{equation}
0 < \mu_\alpha(L_*) = d + \sum_{i=1}^n \beta_i.
\label{3b}
\end{equation}

Since $g_\theta$ is diagonal with respect to this decomposition, $\Phi$ is
either upper or lower diagonal, which means either $L$ or $M$ is
$\Phi$-invariant. But $\alpha$-stability of the pair $(E_*,\Phi)$
implies that
$$\Phi=\left(\begin{array}{cc} 0 & 0 \\ \phi&0 \end{array} \right),$$
where  $0 \neq \phi \in \ParHom(L_*, \widehat{M}_* \otimes K(D)).$
Thus 
$$0 \neq H^0(L_{*}^{\vee} \otimes\widehat{M}_{*} \otimes K(D)) \\
=H^0(L^{\vee}\otimes M\otimes K({\textstyle{\sum}^n_{i=1}} (1-e_i)p_i)).$$
Let $|e| = \sum_{i=1}^n e_i,$ then a necessary condition is that
\begin{equation}
0\leq \deg(L^{\vee}\otimes {M}\otimes K({\textstyle{\sum}^n_{i=1}} (1-e_i)p_i))
=  2(g-1) -2d - |e|. 
\label{3c}
\end{equation}
Now (\ref{3a}) follows from (\ref{3b}) and (\ref{3c}).

We can use
the defining equations for ${\cal E}^0_{\hbox{\scriptsize \it flat}}$ to determine the associated
critical values.
Take $(E_*, \Phi)$ as above, then
$$0 = F_A + [\Phi, \Phi^*] 
 = \left( \begin{array}{cc} F_L - \phi \phi^* & 0 \\ 0 & F_M + \phi^* \! \phi \end{array} \right).$$
Using the Chern-Weil formula for parabolic
bundles (Proposition 2.9 of \cite{biquard}), we get 
$$\mu_\rho(d_A,\Phi) = \frac{1}{4\pi} \| \Phi \|^2 =  \frac{i}{2\pi} \int_X \Tr(\Phi \Phi^*)
= \frac{i}{2\pi} \int_X \phi \phi^* = \frac{i}{2\pi} \int_X F_L = \pardeg (L_*).$$
This completes the proof of (b).

Given $E_* = L_* \oplus M_*$ and $\Phi$ as above, then the zero set of
$\phi$ is a nonnegative divisor of degree 
$$h_{d,e} = \deg (L^{\vee}\otimes {M}\otimes K(\textstyle{\sum^n_{i=1}} (1-e_i)p_i)) =2g-2-2d-|e|$$
on $X,$ which is just an element of ${S}^{h_{d,e}}X.$
Conversely, given a nonnegative divisor of degree $h_{d,e},$
then we obtain a line bundle $U$ of degree $2d+n$ along with a section 
of $U^\vee \otimes K(\textstyle{\sum^n_{i=1}} (1-e_i)p_i))$ vanishing on that divisor.
There are $2^{2g}$ choices of $L$ so that
$U = L^{\otimes 2} \otimes \Lambda^2 E,$
and each choice gives a stable parabolic Higgs bundle $(E_*, \Phi).$
The line subbundle
$L_*$ is canonically determined from $E_*,$ but $\Phi$ is only determined
up to multiplication by a nonzero constant. However, it is easy to see that
$(E_*, \Phi)$ is gauge equivalent to $(E_*, \lambda \Phi)$ for $\lambda \neq 0,$ and (c) now follows.

We now calculate the index $\lambda_{d,e}$ of the critical submanifold ${\cal M}_{d,e},$
which is given by the negative weight space of the
infinitesimal action of $\rho,$ or equivalently, of the gauge transformation $g_\theta.$
Letting $H^0(\ParEnd_0(E)) \cdot \Phi$ be the subspace of Higgs fields of the form
$[\Psi, \Phi]$ for $\Psi \in H^0(\ParEnd_0(E)),$ 
then the subspace
$$W  = H^0(\ParEnd^\wedge_0(E)\otimes K(D)) / H^0(\ParEnd_0(E)) \cdot \Phi$$
is Lagrangian with respect to the complex symplectic form
$$\omega((\xi_1,\phi_1),(\xi_2,\phi_2)) = \int_X \Tr (\phi_2 \xi_1 - \phi_1 \xi_2).$$
So once we determine the weights on $W,$ the weights on the dual space $W^*$ are
given by $1-\nu$ for some weight $\nu$ on $W$
(since $\rho(\theta)^* \omega = e^{i \theta} \omega$).
With respect to the decomposition $E_* = L_* \oplus M_*,$ we
have $$g_\theta = \left( \begin{array}{cc} e^{-i \theta/2} &0 \\ 0 & e^{i\theta/2} \end{array} \right)$$
with weights $(0,1,-1)$ on
$$\ParEnd^{\wedge}_0(E_*) = \ParHom(L_*, \widehat{L}_*) \oplus \ParHom(L_*, \widehat{M}_*) \oplus
\ParHom(M_*, \widehat{L}_*).$$
Further, there are no negative weights on $H^0(\ParEnd_0(E)) \cdot \Phi$ 
and the weights on $W^*$ are $(1,0,2),$ so we get
\begin{eqnarray*}
\lambda_{d,e} &=& 2 h^0(M_*^{\vee}\otimes\widehat{L}_*\otimes K(D)) 
=2(n+2d+g-1+|e|).
\end{eqnarray*}
This completes the proof of (d).
\end{pf}
\subsection{The topology of ${\cal N}^0_\alpha$}
Using the results of the previous section, we deduce the following theorem.
\begin{thm} \label{thm:main}
\begin{itemize}
\item[(a)] If $g>0$ or $g=0$ and $n>3,$ then ${\cal N}^0_\alpha$ is noncompact.
\item[(b)] The Betti numbers of ${\cal N}^0_\alpha$ depend only on the quasi-parabolic structure of $E_*.$
\item[(c)] If $g > 0$ or $g=0$ and $n \geq 3,$ then ${\cal N}^0_\alpha$ is connected and
simply connected.
\end{itemize}
\end{thm}

\begin{pf}
Notice that, whenever $\dim{\cal N}^0_\alpha > 0,$
then for all $(d,e), \; \lambda_{d,e} <  \dim{\cal N}^0_\alpha.$ Thus, the Morse function
$\mu_\rho$ has no maximum value and (a) follows. The only case where
$\dim{\cal N}^0_\alpha = 0$ is, of course, $g=0$ and $n=3.$

We first recall Theorem 3.1 of \cite{bh}.
Let $W= \{\alpha \mid 0 < \alpha_i < \frac{1}{2} \}$ be the
weight space and for any $(d,e),$ define the hyperplane
$H_{d,e} = \{ \alpha \mid d+\beta(\alpha,e) = 0\}.$ 
The set $W \setminus \cup_{d,e} H_{d,e}$ consists of the generic weights, i.e., those for
which stability and semistability coincide.
Suppose $\delta \in H_{d,e},$
then stratifying ${\cal M}^0_{\delta}$ by the Jordan-H\"{o}lder type of the underlying parabolic bundle,
we see that 
$${\cal M}^0_{\delta} = ({\cal M}^0_{\delta} \setminus \Sigma_{\delta}) \cup \Sigma_{\delta},$$
where $\Sigma_{\delta}$ consists of strictly semistable bundles, i.e., semistable
bundles $E_*$ with $\gr E_* = L_* \oplus M_*$ for two parabolic line bundles of parabolic degree zero.
Suppose that $\alpha$ and $\alpha'$
are generic weights
on either side of $H_{d,e}$
and that $\pardeg_\alpha (L_*) < 0.$
If both ${\cal M}^0_{\alpha}$ and  ${\cal M}^0_{\alpha'}$ are nonempty,
then Theorem 3.1 of \cite{bh} states that there are canonical, projective maps
$$\begin{array}{rcl}{\cal M}^0_{\alpha}& & {\cal M}^0_{\alpha'}\\
\phi \!\!\!\! & \searrow \; \swarrow & \!\!\! \phi'\\& {\cal M}^0_{\delta}
\end{array}$$
which are isomorphisms on ${\cal M}^0_{\delta} \setminus \Sigma_{\delta}$ and are
${\Bbb P}^{a}$ and ${\Bbb P}^{a'}$ bundles along $\Sigma_{\delta},$
where $a = h^1(M_*^\vee \otimes L_*)-1$ and $a' = h^1(L_*^\vee \otimes M_*)-1.$
In particular, since $\Sigma_{\delta} = J_X,$ Corollary 3.2 of \cite{bh} gives 
$$P_t({\cal M}^0_{\alpha}) - P_t({\cal M}^0_{\alpha'}) = (P_t({\Bbb P}^{a})-P_t({\Bbb P}^{a'})) P_t(J_X).$$

To prove (b), we must show that $P_t({\cal N}^0_{\alpha}) = P_t({\cal N}^0_{\alpha'})$
for weights on either side of a hyperplane $H_{d,e}.$
Note that $d = \deg L$ and $e = e(L),$ and set $\hat{d}=-n-d$ and $\hat{e}_i = 1-e_i.$
Since 
$$d+\beta(\alpha,e) = \pardeg_\alpha(L) < 0 < \pardeg_{\alpha'}(L) = d + \beta(\alpha',e),$$
and 
$\hat{d} +\beta(\alpha', \hat{e}) < 0 < \hat{d}+\beta(\alpha,\hat{e}),$
it follows that the indexing sets of $(d,e)$ satisfying (\ref{3a})
for ${\cal N}^0_\alpha$ and ${\cal N}^0_{\alpha'}$ are identical except 
for $(d,e)$ and $(\hat{d},\hat{e})$ listed above; the pair $(d,e)$ satisfies (\ref{3a}) for $\alpha$ but
not for $\alpha'$ and vice versa for $(\hat{d},\hat{e}).$
Thus, we claim
\begin{eqnarray*}
0 
&=& P_t({\cal M}^0_{\alpha}) - P_t({\cal M}^0_{\alpha'}) 
+ t^{\lambda_{d,e}} P_t({\cal M}_{d,e})- t^{\lambda_{\hat{d},\hat{e}}} P_t({\cal M}_{\hat{d},\hat{e}}), 
\end{eqnarray*}
which, setting $\Delta = t^{\lambda_{\hat{d},\hat{e}}} P_t({\cal M}_{\hat{d},\hat{e}}) - t^{\lambda_{d,e}} P_t({\cal M}_{d,e})$  is equivalent to 
\begin{equation} \label{3k}
\Delta = \frac{( t^{2a'+2} -t^{2a+2}) (1+t)^{2g}}{1-t^2}.
\end{equation}
First, we compute
$$h_{d,e} = 2g-2-2d -|e|, \quad  \lambda_{d,e} = 2(n+2d+g-1+|e|), $$
$$h_{\hat{d},\hat{e}} = 2g-2+n+2d+|e|, \quad  \lambda_{\hat{d},\hat{e}} = 2(g-1-2d-|e|).$$
Next, notice that if $h> 2g-2,$  then
$P_t(\widetilde{S}^h(X))= P_t(S^h(X))$ (see p. 98 of \cite{hitchin}).
But both $h_{d,e}$ and $h_{\hat{d},\hat{e}}$ are greater than $2g-2,$ 
which we see as follows.
Since $\frac{e_i}{2} \leq \beta_i(\alpha,e) \leq \frac{1+e_i}{2},$
we have $\frac{|e|}{2} \leq \sum_{i=1}^n \beta_i(\alpha,e) \leq \frac{n+|e|}{2}.$
It now follows that
$ 2d + |e| < 2d + 2\beta(\alpha,e)  < 0$ and $2d+n+|e| > 2 d + 2 \sum_{i=1}^n \beta(\alpha',e) > 0.$

Now use the result of \cite{macdonald} to interpret
$P_t(S^h X)$ as the coefficient of $x^h$ in $$\frac{(1+xt)^{2g}}{(1-x)(1-xt^2)},$$ 
and compute in terms of residues to  see 
\begin{eqnarray*}
\Delta&=& t^{\lambda_{\hat{d},\hat{e}}}P_t(S^{h_{\hat{d},\hat{e}} }X) - t^{\lambda_{d,e}}P_t(S^{h_{d,e}}X) \\
&=& \Res_{x=0} \left( \frac{t^{\lambda_{\hat{d},\hat{e}}}}{x^{h_{\hat{d},\hat{e}}+1}} 
- \frac{t^{\lambda_{d,e}}}{x^{h_{d,e}+1}}
\right)\left( \frac{(1+xt)^{2g}}{(1-x)(1-xt^2)}\right).
\end{eqnarray*}
This last function is analytic at $x=\infty$ and has a removable singularity at $x=1/t^2,$ thus
\begin{eqnarray*} 
\Delta &=& -\Res_{x=1} \left( \frac{t^{\lambda_{\hat{d},\hat{e}}}}{x^{h_{\hat{d},\hat{e}}+1}}
- \frac{t^{\lambda_{d,e}}}{x^{h_{d,e}+1}}
\right)\left( \frac{(1+xt)^{2g}}{(1-x)(1-xt^2)}\right)\\
&=& \frac{(t^{\lambda_{\hat{d},\hat{e}}} - t^{\lambda_{d,e}})(1+t)^{2g}}{1-t^2}.
\end{eqnarray*}
But we can compute directly that
$2a'+2 = \lambda_{\hat{d},\hat{e}}$ and that $2a+2 = \lambda_{d,e}$ and (\ref{3k}) follows. 
This proves (b) in case both ${\cal M}^0_\alpha$ and ${\cal M}^0_{\alpha'}$ are nonempty.
In case one of the moduli is empty, we use the following lemma (see the remark).

To prove (c), we use the fact that
${\cal M}^0_\alpha$ is connected and simply-connected, which follows 
for $g=0$ from \cite{bauer} and for $g \geq 1$ from \cite{boden1}.
Since $\lambda_{d,e}$ is always even,
(c) will follow if $\lambda_{d,e} >0$ for all $(d,e).$
This is true if ${\cal M}^0_\alpha \neq \emptyset.$
However, if $g=0$ we must be careful since there are weights
$\alpha$ with ${\cal M}_\alpha = \emptyset.$
In that case, we must show that there is a unique pair $(d,e)$ with
$\lambda_{d,e} = 0,$ and also that ${\cal M}_{d,e}$ is connected and simply connected.
This is the content of the following lemma. 
\end{pf}
\begin{lem}
\begin{enumerate} 
\item[(i)] If $g \geq 1,$ then $\lambda_{d,e} > 0$ for every $(d,e)$ satisfying \rom(\ref{3a}\rom). 
\item[(ii)] If $g = 0$ and $n \geq 3,$ then there is at most one pair $(d,e)$ satisfying
\rom(\ref{3a}\rom) with $\lambda_{d,e} = 0.$ Such a pair $(d,e)$ exists 
if and only if ${\cal M}_\alpha = \emptyset,$
and in that case, ${\cal M}_{d,e} = {\Bbb P}^{n-3}.$ Here, $ {\cal M} = {\cal M}^0$ since $g=0.$
\end{enumerate}
\end{lem}
{\it Remark. \,} 
We now explain why this lemma proves part (b) of the Proposition when
one of the moduli is empty.  Suppose ${\cal M}_{\alpha}=\emptyset,$ 
then it follows that the moment map $\mu_\rho$ is positive with minimum value
$d+\sum_{i=1}^n \beta(\aa,e)$ for the pair $(d,e)$ identified in part (ii) of the lemma.
Since $(d,e)$ does not satisfy (\ref{3a})
for $\alpha',$ 
$H_{d,e}$ is the relevant hyperplane.
This identifies the birth and death strata as ${\cal M}_{\aa'}$ and ${\cal M}_{d,e},$
and thus all the other strata for $\aa$ and $\aa'$ are identical.
The rest follows from the fact that
${\cal M}_{\alpha'} = {\Bbb P}^{n-3},$ 
first proved by Bauer \cite{bauer}.
\begin{pf}
Suppose that $\lambda_{d,e} =0$ for a pair $(d,e)$ satisfying (\ref{3a}).
We first show that $g=0.$ Recall that $\beta_i(\alpha,e) = e_i +(-1)^{e_i} \alpha_i.$
Using the fact that $0 = \lambda_{d,e} = n+2d+g+|e|-1,$ 
the condition (\ref{3a}) and the inequality $\beta_i(\alpha,e) < \frac{e_i +1}{2},$ we see that
\begin{equation} \label{3d}
\frac{n+|e|+g-1}{2} < \sum_{i=1}^n \beta_i(\alpha, e) <  \frac{n+|e|}{2}.
\end{equation}
This is only possible if $g=0,$ which we now assume.

Setting $\gamma_i = 1- \beta_i = (1-e_i)(1-\alpha_i)+e_i \alpha_i,$ then
equation (\ref{3d}) is equivalent to
$$\frac{n-|e|}{2} < \sum_{i=1}^n \gamma_i < \frac{n-|e|+1}{2}.$$
Writing $\gamma_i = \frac{1-e_i}{2} + (1-e_i)(\frac{1}{2}-\alpha_i) + e_i \alpha_i,$
we get immediately 
\begin{equation} \label{3e}
0< \sum_{i=1}^n (1-e_i)(\frac{1}{2}-\alpha_i) + e_i \alpha_i < \frac{1}{2}. 
\end{equation}
The advantage of the (\ref{3e}) is that each summand is positive.

We now prove uniqueness of the pair $(d,e).$
If $\lambda_{d',e'} = 0$ for $(d',e') \neq (d,e),$
then it follows that $|e|-|e'| = 2(d'-d)$ is even, which implies that $e_i \neq e'_i$
for at least two $i,$ which we assume (wlog) to include $i=1,2.$
Now $(\alpha,e)$ and $(\alpha, e')$ both satisfy the inequality (\ref{3e}).
Add them together and notice that
since $e_1 \neq e'_1$ and $e_2 \neq e'_2,$
the sum of the left hand sides is at least
$ \alpha_1 + (1/2-\alpha_1) + \alpha_2 + (1/2-\alpha_2) = 1,$
which violates the (summed) inequality and therefore gives a contradiction.

It follows from $\lambda_{d,e}=0$ and $g=0$ that $n+|e|-1$ is even and $h_{d,e} = n-3.$ 
Thus ${\cal M}_{d,e} = S^h X =S^h {\Bbb P}^1= {\Bbb P}^{n-3}.$ 
The rest of the lemma follows from the the inequality (\ref{3d}), together with the
following proposition,
which we have chosen to state as it is of independent interest.
\end{pf}
\begin{prop} \label{prop:null-ch}
If $g=0,$ then the moduli space ${\cal M}_\alpha \neq \emptyset \Leftrightarrow$ 
\begin{equation}\label{3f}
\sum_{i=1}^n e_i + (-1)^{e_i}\alpha_i < \frac{n+|e|-1}{2}.
\end{equation}
for every $e=(e_1,\ldots,e_n), \; e_i \in \{0,1\},$ with $n-|e|+1$ even.
\end{prop}
\noindent
{\it Remark. \,} For $n=3,$ ${\cal M}_\alpha$ is either empty or a point.
In this case, the proposition can be verified directly by
comparing the inequalities (\ref{3f})
to the well-known fusion rules (or the quantum Clebsch-Gordan conditions):
$${\cal M}_\alpha \neq \emptyset  \Leftrightarrow
|\alpha_1 - \alpha_2| \leq \alpha_3 \leq \min(\alpha_1 + \alpha_2, 1-\alpha_1-\alpha_2).$$
\begin{pf}
Like the proof of part (b) of the theorem, we shall use the techniques of \cite{bh}.
Recall the weight space $W = \{ \alpha \mid 0 \leq \alpha_i \leq 1/2 \}$
and the hyperplanes $H_{d,e} = \{ \alpha \mid d+\beta(\alpha,e)=0 \}$ defined earlier.
We call connected components of $W \setminus \cup_{d,e} H_{d,e}$ {\it chambers}.
A chamber $C$ is called {\it null} if the associated moduli space ${\cal M}_\alpha$ is empty in genus 0
for every $\alpha \in C.$
The proposition follows once we show that every null chamber is given by
$C_{d,e} = \{ \alpha \mid d+\beta(\alpha,e)>0 \},$ where $2d=1-n-|e|.$ 

Associated to the configuration of hyperplanes in $W$ is a graph
with one vertex for each chamber and an edge between two vertices whenever the two chambers
are separated  by a  hyperplane.
We shall see that in terms of this graph, null chambers have valency one.
The (unique) hyperplane separating a null chamber from the rest of $W$ is
called a {\it vanishing wall.}
If $\delta \in H_{d,e},$ a vanishing wall, and 
$\alpha, \alpha'$ are nearby weights on either side of $H_{d,e},$
then the proof of Proposition 5.1 of \cite{bh} shows that 
${\cal M}_\delta = \Sigma_\delta$ and, assuming that ${\cal M}_{\alpha'} = \emptyset,$
the map $\phi$ is a fibration with fiber
${\Bbb P}^{a},$ where $a = h^1(M_*^\vee \otimes L_*)-1.$
Moreover, $h^1(L_*^\vee \otimes M_*)=0$ and this last equation in fact characterizes vanishing walls.

We claim that every vanishing hyperplane is given by $H_{d,e}$ for $2d=1-n-|e|.$
For if $d=\deg L$ and $e=e(L),$ then direct computation shows that
$h^1(L_*^\vee \otimes M_*) =2d+n+|e|-1.$ 
On the other hand, if 
$n+|e|-1$ is even and $d = \frac{1-n-|e|}{2},$
then $H_{d,e}$ is a vanishing hyperplane. 

Along $H_{d,e},$ the relevant line bundles of parabolic degree 0
are given by $L_* = {\cal O}_X (\frac{-n-|e|+1}{2})[-\beta]_*$
and $M_* = {\cal O}_X (\frac{-n+|e|-1}{2})[-\gamma]_*,$
where $\delta \in H_{d,e}, \beta= \beta(\delta,e)$ and $\gamma_i = 1- \beta_i.$ 
Since 
$h^1(L_*^\vee \otimes M_*)=0$ and $h^1(M_*^\vee \otimes L_*)=n-2,$
it follows that the null chamber is defined by
$C_{d,e}=\{ \alpha \mid \beta(\alpha,e) > \frac{n+|e|-1}{2} \}.$
To verify that this is indeed a chamber, we prove that no other hyperplane cuts through
$C_{d,e}.$ This will also show 
that null chambers have valency one in
the graph associated to the configuration of hyperplanes.

Suppose to the contrary that $\alpha \in H_{d',e'} \cap C_{d,e}.$ 
Then we have
$\sum (-1)^{e_i} \alpha_i > \frac{n-|e|-1}{2}$
and $\sum (-1)^{e'_i} \alpha_i = -|e'| -d' =k \in {\Bbb Z}.$
If $e_i = e_i' = 0,$
then $((-1)^{e_i} + (-1)^{e_i'}) \alpha_i < 1$ and in all other cases,
$((-1)^{e_i} + (-1)^{e_i'}) \alpha_i  \leq 0.$ Using a similar property for $e''= 1 -e',$ we see
\begin{eqnarray*}
\frac{n-|e|-1}{2} + k &<& \sum_{i=1}^n((-1)^{e_i} + (-1)^{e_i'}) \alpha_i < \sum_{e_i = e'_i=0} 1,\\
\frac{n-|e|-1}{2} -k &<& \sum_{i=1}^n((-1)^{e_i} + (-1)^{e_i''}) \alpha_i < \sum_{e_i = e''_i=0} 1.
\end{eqnarray*}
These are strict inequalities of integers,
so after adding one to the left hand sides and summing the two inequalities (which are no longer strict),
we see 
$n-|e|+1 \leq \sum_{e_i=0} 1 = n-|e|,$ a contradiction.
\end{pf}

\subsection{The Betti numbers of the moduli space of parabolic Higgs bundles} 
The results of the previous section show that the Betti numbers of ${\cal N}^0_\alpha$
depend only on the genus $g$ and number $n$ of parabolic points.
In this section, we 
give a formula for the Poincar\'e polynomial of ${\cal N}^0_\alpha.$ 
Such a general calculation is not possible for $P_t({\cal M}^0_\alpha)$
without first specifying $\alpha,$ 
so take $\alpha=(\frac{1}{3},\ldots, \frac{1}{3^n}).$ Using Proposition \ref{prop:null-ch}
(taking $e=(0,1,\ldots,1)$) it is clear that $\alpha$ lies in a null chamber.
We could calculate $P_t({\cal M}^0_\alpha)$ using the Atiyah-Bott procedure for
parabolic bundles as in \cite{boden1}, but there is an easier method which
exploits the fact that $\alpha$ lies in a null chamber.
First of all, using the results of \S 6.4 in \cite{boden1}, we get
$$P_t({\cal M}^0_\alpha) = \frac{(1+t^2)^{n-1}(1+t^3)^{2g}}{(1-t^2)^2}- \frac{(1+t)^{2g}}{(1-t^2)} \sum_{\lambda,e}t^{2 d_{\lambda,e}}.$$
Note that $d_{\lambda,e}$ depends on $g$
($d_{\lambda,e} =d_{\lambda,e}(g=0)+g$), but the indexing set $\{ \lambda,e\}$ is independent of $g$.
Since ${\cal M}^0_\alpha(g=0) = \emptyset,$ this determines the sum and we see that 
$$P_t({\cal M}^0_\alpha) = \frac{(1+t^2)^{n-1}\left( (1+t^3)^{2g}-t^{2g}(1+t)^{2g}}{(1-t^2)^2\right)}.$$
It follows from Theorem \ref{thm:morse} that
$$P_t({\cal N}^0_\alpha) = P_t({\cal M}^0_\alpha) + \sum_{d,e} t^{\lambda_{d,e}} P_t({\cal M}_{d,e}),$$
where the sum is taken over $(d,e)$ satisfying (\ref{3a}),
which,
for our choice of $\alpha,$ is simply
$e_1 -|e| \leq d \leq [g-1 -\frac{|e|}{2}]$, where $[x]$ is the greatest integer less than $x.$
Setting $j=2d+n+|e|-1,$ then $j$ satisfies:
$$n+2e_1 -|e|-1 \leq j \leq 2g+n-3 \quad \hbox{ and } \quad j-n-|e|+1 \hbox{ is even.}$$
Also $\lambda_{d,e} = 2(g+j)$ and $h_{d,e} = 2g+n-j-3.$

Fixing $e_1$ and $|e|,$ for each $d,$ there are $\left({{n-1}\atop{|e|-e_1}}\right)$ strata
given by the choice of $e.$ Thus, for each $j,$ there are 
$q_j = \sum_{i=0}^{j} \left({{n-1}\atop{i}}\right)$ strata
(note that $q_j = 2^{n-1}$ for $j \geq n-1$) and we see
\begin{eqnarray*}
\sum_{d,e} t^{\lambda_{d,e}} P_t({\cal M}_{d,e}) \!\!&=& \!\!
\sum_{|e|=0}^n \left({{n-1}\atop{|e|-e_1}}\right) 
\sum_{d=e_1 - |e|}^{[g-1-|e|/2]} t^{\lambda_{d,e}} P_t(\widetilde{S}^{h_{d,e}}X) \\
\!\! &=& \!\! \sum_{j=0}^{2g+n-3} q_j t^{2(g+j)} P_t(\widetilde{S}^{2g+n-j-3}X)\\
\!\! &=& \!\! \sum_{j=0}^{n-2} q_j t^{2(g+j)} P_t(\widetilde{S}^{2g+n-j-3}X)
+ \sum_{j=0}^{2g-2} 2^{n-1}t^{2(g+n+j-1)} P_t(\widetilde{S}^{2g-j-2}X). 
\end{eqnarray*}
We refer to the last two sums by $\widetilde{S}_1$ and $\widetilde{S}_2.$
Using the Binomial Theorem and the general formula (p. 98 of \cite{hitchin}) 
$\P_t(\widetilde{S}^{h}X) = (2^{2g}-1)\left({{2g-2}\atop{h}}\right)t^h +\P_t({S^{h}X}),$
we see that
\begin{eqnarray*}
\widetilde{S}_1 &=& \sum_{j=0}^{n-2} q_j t^{2(g+j)} P_t({S}^{2g+n-j-3}X) = S_1,\\
\widetilde{S}_2 &=& \sum_{j=0}^{2g-2} 2^{n-1}t^{2(g+n+j-1)} P_t({S}^{2g-j-2}X) 
            + \sum_{j=0}^{2g-2} 2^{n-1} (2^{2g}-1) \left({{2g-2}\atop{j}} \right) t^{4g+2n+j-4}\\
&=& S_2 +  2^{n-1} (2^{2g}-1) t^{2(2g+n-2)}(1+t)^{2g-2},
\end{eqnarray*}
where $S_1$ and $S_2$ are the sums obtained by removing the tildes from the 
summands of $\widetilde{S}_1$ and $\widetilde{S}_2.$
According to a result of \cite{macdonald}, 
$\P_t({S^{h}X})$ is the coefficient of $x^h$ in $$\frac{(1+xt)^{2g}}{(1-x)(1-xt^2)}.$$
This allows us to evaluate $S_i$ as follows:
\begin{eqnarray*}
S_1 &=& \Res_{x=0} 
\left(\sum_{j=0}^{n-2}\frac{q_j t^{2(g+j)} (1+xt)^{2g}}{x^{2g+n-j-2}(1-x)(1-xt^2)}\right),\\
S_2 &=& \Res_{x=0}  
\left( \frac{2^{n-1}t^{2(g+n-1)}(1+xt)^{2g}}{x^{2g-1}(1-x)(1-xt^2)^2}\right).
\end{eqnarray*}  
But each of these rational functions is analytic at $x = \infty,$
so we can use the Cauchy Residue Formula to evaluate instead at the poles $x=1$ and $x=1/t^2.$
Letting $Q_n(t) = \sum_{k=0}^{n-2} q_k t^{2k}$ and noticing that
$Q_n(1) = \sum_{k=0}^{n-2} q_k = 2^{n-2}(n-1),$
we get
\begin{eqnarray*}
S_1 &=&  \left(Q_n(t) t^{2g} - 2^{n-2}(n-1) t^{2(2g+n-2)}\right) \frac{(1+t)^{2g}}{(1-t^2)},\\
S_2 &=&  2^{n-1}\left(t^{2(g+n-1)} + t^{4g+2n-3} \left( (2g-1)t - 2g \right) \right) \frac{(1+t)^{2g}}{(1-t^2)^2}.
\end{eqnarray*}
But since $Q_n(t)(1-t^2) + 2^{n-1}t^{2(n-1)} = (1+t^2)^{n-1},$
it follows that
\begin{eqnarray*} 
P_t({\cal N}^0_\alpha) &=& P_t({\cal M}^0_\alpha) + \widetilde{S}_1 + \widetilde{S}_2 \\
&=& P_t({\cal M}^0_\alpha) + S_1 + S_2 + 2^{n-1}(2^{2g}-1) t^{2(2g+n-2)}(1+t)^{2g-2}\\\\
&=& \frac{(1+t^3)^{2g}(1+t^2)^{n-1} + 2^{n-1}t^{2n+4g-3}(1+t)^{2g}[(2g-1)t - 2g]}{(1-t^2)^2}\\
& & -\frac{2^{n-2}(n-1)t^{2n+4g-4}(1+t)^{2g}}{1-t^2} + 2^{n-1}(2^{2g}-1) t^{4g+2n-4}(1+t)^{2g-2}.
\end{eqnarray*}
Evaluating this at $t=-1$ 
shows that the Euler characteristic of ${\cal N}^0_\alpha$ is given by
$$\chi({\cal N}^0_\alpha)
= \begin{cases} (n-1)(n-2)2^{n-4} & \text{if } g=0, \\ 3\cdot2^n & \text{if } g=1, \\ 0 & \text{if } g \geq 2
\end{cases}$$

Theorem \ref{thm:main}
would lead one to believe that the diffeomorphism type of
${\cal N}_\alpha^0$  depends only on the quasi-parabolic structure.
We conjecture this is true in general.
Subsequent to the writing of this paper,
this conjecture was proved by H. Nakajima in rank
two \cite{nakajima}.

\end{document}